# Charge-Discharge Coupling Strategy for Dispatching Problems with Electric Tractors at Airports

Danwen Bao[1], Ziqian Zhang[2], Di Kang[3*]

Corresponding Author: Di Kang

[1]College of Civil Aviation, Nanjing University of Aeronautics and Astronautics, Nanjing, 211106, China

E-mail: baodanwen@nuaa.edu.cn

[2]College of Civil Aviation, Nanjing University of Aeronautics and Astronautics, Nanjing, 211106, China

E-mail: 2514323741@qq.com

[3]Department of Civil and Environmental Engineering, Rutgers University, Piscataway, NJ, 08854, USA

E-mail: di.kang@rutgers.edu

A B S T R A C T

Airports worldwide are actively promoting the transition of ground service vehicles from traditional fuel-powered vehicles to electric vehicles. The key to the successful implementation of this transition lies in the development of efficient electric vehicle dispatching models that comprehensively consider the charge-discharge processes of electric vehicles. However, due to the nonlinear characteristics of charge-discharge processes, finding precise solutions poses a significant challenge. Previous researchers have often used traditional energy consumption models and constant charging rates to simplify calculations, but this has resulted in inaccurate estimates of the remaining battery charge level. Furthermore, the lack of diverse pacing and charging strategies for airport ground service vehicles necessitates more adaptable solutions to enhance operational efficiency. To address these challenges, this paper uses airport electric tractors as a case study, develops an accurate model that takes into account the start-stop process and a piecewise linear charging function, designs an improved genetic algorithm that incorporates a greedy algorithm and an adaptive strategy, and develops charge-discharge coupling strategies for different configuration scenarios at Nanjing Lukou Airport to meet current and future needs. The research results indicate that compared to traditional genetic algorithms, the proposed improved genetic algorithm significantly enhances solution accuracy and convergence speed. Additionally, with the increase in flight scale, airports can appropriately enhance their charging strategies; airports with dispersed aircraft stands should devise higher pacing strategies compared to those with dense aircraft stands.

## 1. Introduction

In recent years, the issue of global warming caused by large-scale greenhouse gas emissions has become increasingly severe, and the problem of reducing carbon dioxide emissions has gathered widespread international attention. The aviation industry, as a significant component of global economic activities, is facing an increasingly urgent need to reduce carbon emissions. According to historical data from the International Energy Agency, before the outbreak of the COVID-19 pandemic in 2019, the global aviation industry produced 918 million tons of carbon emissions, accounting for 2% of the total global carbon emissions. Traditional airport ground service vehicles are the second-largest source of carbon emissions in the aviation industry, trailing only behind emissions from the aircraft themselves. These vehicles, such as baggage conveyors, shuttle buses, tractors, and others, serve various functions at airports and consume a substantial amount of fuel. Converting these vehicles from traditional fossil fuel power to electric power could prevent serious environmental pollution.

In the contemporary era, research on this transition primarily focuses on the public transportation and logistics industry. Depending on the study approaches, it can be categorized into three main categories: vehicle routing problem for electric vehicles (Li et al., 2018; Basso et al., 2019; Yao et al., 2020; Tang et al., 2023), site selection for charging stations (He, Yin, and Zhou, 2015; Tu et al., 2016; Wu and Sioshansi, 2017; Tadayon-Roody, Ramezani, and Falaghi, 2021),and intelligent charging strategies (Hiermann et al., 2016; Amini, Kargarian, and Karabasoglu, 2016; Khalkhali and Hosseinian, 2020; Zweistra, Janssen, and Geerts, 2020). One of the hot research topics is the trajectory optimization of electric vehicles, but these studies generally have the following caveats: firstly, numerous studies make the assumption that the battery discharge model is linear, which does not align with the actual nonlinear discharge behavior. Similarly, linear charging models are often used in research, while the charging rate in reality is also nonlinear. This mismatch may lead to charging plan failures, vehicles being unable to complete tasks, and consequences such as battery degradation due to undercharging storage.

To be more specific, previous studies often assumed that electric vehicles maintain a constant speed throughout their journeys (Bao et al., 2023), without accounting for acceleration or deceleration. This overlook of the start-stop processes has led to inaccuracies in the discharge models of electric vehicles. In practice, electric vehicles consume more energy during the initial acceleration phase and recuperate a certain amount of energy during braking (Basso et al., 2019). In contrast to regular road traffic, airport service vehicles tend to have shorter driving trajectories and experience more frequent instances of starting and stopping. Consequently, neglecting start-stop processes in the discharging model analysis could lead to a notable discrepancy between the calculated and real energy consumption values. Moreover, in contrast to standard electric vehicles, airport electric vehicles frequently supply specialized functions for aircraft operations. For instance, Commercial Aircraft Corporation of China, Ltd. has developed an electric tractor, CN214986146U, which not only provides towing services for aircraft but also supplies power to the aircraft's lighting and air conditioning systems (van Baaren, 2019). These vehicles must distinguish between energy consumption during driving and energy consumption during service operations, rendering conventional energy consumption models obsolete.

On the other hand, there are also issues with the charging models for electric vehicles. Previous research often assumed that charging rates are constant (Hof, Schneider, and Goeke, 2017) or simply linearly related to the charging amount (Lin, Ghaddar, and Nathwani, 2022). Indeed, due to variations in terminal voltage and current throughout the charging process, the charging function ought to be nonlinear. An inaccurate fitting of the charging function can result in incorrect estimates of charging durations, ultimately leading to impractical dispatching strategies. In other transportation fields, charging strategies exhibit flexibility and variability (Zhang et al., 2018), but the airport commonly employs a full charging strategy. This practice stems from the earlier assumption that charging amount is linearly related to charging time, where achieving 100% of the maximum charging level was considered the most optimal level. However, when considering a nonlinear charging function, the optimal charging level is not necessarily 100%. Based on an actual charging curve depicting the relationship between battery charge level and charging time, the charging rate remains relatively high until reaching a threshold, beyond which the charging rate decreases significantly (Montoya et al., 2017). Therefore, a trade-off emerges between the battery charge level and the frequency of necessary charging sessions: charging the battery until it reaches a low level and stopping at that threshold, despite the fast-charging rate, results in frequent recharges for the vehicles. On the other hand, if we mandate that vehicles charge to a high battery charge level, it does reduce the number of required charging sessions. However, this approach also results in slower speeds during the latter part of the charging process, which ultimately leads to excessively long individual charging durations. Thus, for airports with varying throughputs, the best charging strategy is no longer a one-size-fits-all approach, and the full charging strategy used in previous work does not apply universally to all situations.

The preceding paragraphs underscore the existing issues with energy consumption and charging models, leading to impractical electric vehicle dispatching plans and diminishing the applicability of prior research findings. Airport electric vehicles present substantial differences when compared to conventional and general-purpose electric vehicles. In comparison to general-purpose electric vehicles, airport electric vehicles are dedicated to ground service tasks and thus have different battery capacity and operational characteristics (Brevoord, 2021). Therefore, the electric vehicle dispatching models developed for general-purpose electric vehicles cannot be readily applied in the airport. This necessitates the development of a specialized dispatching model that considers the constraints associated with the ground service processes of airport electric vehicles. Within the mathematical model, variables like charging timing, maximum charging level, maximum driving speed, and vehicle servicing arrangements need to be determined to ensure efficient dispatch without unwarranted energy depletion.

The dispatching model developed in this paper advances prior approaches in the following aspects. Firstly, this model considers the nonlinearity in energy consumption caused by vehicle start-stop processes and service processes, addressing the errors arising from the linear assumption. Furthermore, an actual charging dataset is utilized to derive a precise and accurate charging curve, which served as the foundation for introducing the piecewise linear charging function. In contrast to previous research, this marks the first instance of applying the piecewise linear charging function to electric vehicles at airports. Finally, this study introduces, for the first time, the combined consideration of both discharging and charging processes for airport electric vehicles, presenting a mathematical dispatching model with a charge-discharge coupling strategy. This approach proves advantageous in addressing the limitation of having separate vehicle pacing and charging strategies, as well as providing differentiated strategies for airports of different configurations and throughputs.

The main contributions of this paper are as follows:

- This paper implements a simple temporal network approach. It considers sequential constraints and time limitations in airport operations, enabling the establishment of precise time windows for aircraft towing services, which play a vital role in optimizing real-world airport operations.
- It develops a sophisticated energy consumption model that accounts for electric tractor start-stop processes and service status. Additionally, this paper establishes a more realistic airport electric tractor dispatching model by leveraging a piecewise linear charging function.
- This paper introduces a novel concept of charge-discharge coupling strategy, which not only reduces costs but also enhances vehicle utilization rates. This innovation offers theoretical guidance for optimizing the dispatching and charging strategies of electric tractors across airports.
- It devises an improved genetic algorithm that combines greedy principles with adaptive strategies, effectively addressing the challenges posed by excessive model constraints (resulting in poor initial solution quality and susceptibility to local optima). This improvement results in enhanced solution accuracy and faster convergence of the algorithm.

The remaining sections of this paper are organized as follows: Section 2 reviews relevant literature. Section 3 introduces the new methodology adopted in this paper and presents a novel concept of charge-discharge coupling strategy. Section 4 establishes an airport electric tractor dispatching model with time windows and a nested complex energy consumption model. Section 5 outlines the solution algorithm designed for the model. Section 6 conducts a sensitivity analysis of the algorithms and models and presents case study results. Section 7 tailors the present and future charge-discharge coupling strategies for the two scenarios of Nanjing Lukou Airport. Section 8 provides a comprehensive summary of the entire paper.

## 2. Literature review

This section first reviews the literature on the energy consumption of electric vehicles, then lists several charging behaviors, and finally describes the current state of research on airport ground service vehicle dispatching.

In terms of electric vehicle energy consumption, Bektas et al. (2011) first applied an energy consumption model to a vehicle routing problem(VRP), taking into account the effects of load, speed, and road gradient on fuel consumption. Goeke et al. (2015) considered the efficiency of the conversion between mechanical energy and electrical energy, and for the first time applied an energy consumption model to an electric vehicle routing problem(EVRP). Turkensteen (2017) investigated the deviation of the predicted energy consumption from the actual energy consumption under the assumption of constant speed and concluded that this deviation could be up to double. Basso et al. (2019) summarized air resistance, rolling resistance, power system losses and auxiliary system losses as the four major influences on the energy consumption of electric vehicles, as well as calculating an energy consumption integral formula considering the start-stop process. Most of the existing studies have established simple linear energy consumption models, and only a few have considered the effects of load, speed, acceleration and deceleration of electric vehicles. In real life, load positively correlates with energy consumption of electric vehicles. Electric vehicles consume more energy in the start-up phase in order to change the driving state, and a certain amount of energy is recovered in the braking phase due to the motor reversal. The energy consumption of electric vehicles in the travel and service states also differs. The large deviation between the predicted and actual values of energy consumption can lead to the failure of charging schedules, the inability of vehicles to fulfil their tasks, and battery damage from undercharged storage. This is unacceptable for airports that are highly motivated by low costs and low delays.

In terms of charging behaviors of electric vehicles, there are two main categories: full charging and partial charging. In full charging, the vehicle is charged halfway to a fully charged state. Full charging is further subdivided into two categories according to the charging time assumption: scholars represented by Conrad et al. (2011), Adler et al. (2014), Hof et al. (2017) assumed a constant charging time; other scholars represented by Goeke et al. (2015), Hiermann et al. (2016), Lin et al. (2022) assumed that the charging time is a linear function of the battery charge level. In partial charging, the battery level of vehicles is the decision variable and does not have to be charged to full. Almost all scholars represented by Felipe et al. (2014), Montoya et al. (2017), Wang et al. (2013) assumed that the charging time is a linear approximate function of the battery charge level. Most of the existing studies assume the charging time to be constant or a simple linear correlation with the battery charge level. However, in general, the charging function should be nonlinear due to the changes in terminal voltage and current during the charging process. In this regard, most scholars have used optimistic or pessimistic criteria to fit the linear part of the nonlinear charging function, and only a few scholars have considered the use of piecewise linear approximation methods. Besides this, the field of airport vehicle routing has not yet been solved by using a more realistic piecewise linear charging function model.

In terms of airport ground service vehicle dispatching, there are two main categories: fuel vehicle dispatching and electric vehicle dispatching. In fuel vehicle dispatching, Du (2014), Guo et al. (2020), Han et al. (2023) and other scholars have studied the single-type ground service vehicle dispatching problem. All of them transformed the problem into a vehicle routing problem with a time window, and solved it by using different objective functions and algorithms: From the viewpoint of objective function, there are single-objective and multi-objective problems, which include the shortest total distance travelled, the least operating cost, the least delay, etc. From the perspective of solution algorithms, there are exact algorithms and heuristic algorithms. Scholars such as Zhao et al. (2019) and Zhu et al. (2022) have studied the multi-type ground service vehicle dispatching problem, which mainly involves the collaborative dispatching of two kinds of vehicles, and the setting of the expected start moment and duration of ground service operations based on manual experience values in static scenarios. In electric vehicle dispatching, only a very few scholars, such as Brevoord (2021), have combined the ground service process and electric characteristic constraints to establish an airport electric ground service vehicle dispatching model with the objective of the shortest charging time. Most of the studies focus on the urban public transport field and logistics field, represented by Wen et al. (2016), Ma et al. (2018) and Lech et al. (2022), which establish a large number of constraints and improve the traditional heuristic algorithms by adopting optimization operators and novel coding. There are a lot of research papers on airport fuel ground service vehicle dispatching, but the constraints are numerous and difficult to solve, which makes it hard to satisfy the characteristics of high real-time requirement of airport vehicle dispatching plans. In addition, there is a lack of research on the dispatching of electric airport ground service vehicles. Electrified airport ground service vehicles are different from typical electric vehicles, and most of them have their own ground service tasks. There are also variations between the two with respect to battery capacity, charging pile location, and travelling routes.

## 3. Theoretical background

This section consists of the following subsections: (1) the impact of start-stop processes on energy consumption during the discharging process, (2) the impact of using a piecewise linear charging function on charging time during the charging process, and (3) based on findings from (1) and (2), we ultimately propose a charge-discharge coupling strategy.

*3.1. Energy consumption during start-stop processes*

Traditional energy consumption models overlook the start-stop processes of electric vehicles, assuming that the cumulative energy consumption is simply linearly related to the distance traveled. To standardize the unit of energy consumption, we introduce a constant factor, $\epsilon = \frac{1}{3{,}600{,}000}$, representing the conversion difference between joules and kilowatt-hours (kWh). The equation for the cumulative energy consumption within a hundred meters, denoted as $q_{cum}$, is shown as equation (1):

$$q_{cum} = \left(C_r \cdot m \cdot g + \frac{1}{2} \cdot \rho \cdot A \cdot C_d \cdot v^2\right) \cdot l \cdot \emptyset^d \cdot \varphi^d \cdot \epsilon, 0 \leq \ l \leq 100 \tag{1}$$

Where:

$C_r$: the rolling friction coefficient.

$m$: the mass of the vehicle (kg).

$g$: the acceleration due to gravity (N/kg).

$\rho$: air density (kg/m$^3$).

$A$: the frontal area of the vehicle (m$^2$).

$C_d$: the air resistance coefficient.

$v$: the speed of the vehicle (m/s).

$l$: the distance traveled by the vehicle (m).

$\emptyset^d$: the output efficiency of the drive motor when releasing electrical energy.

$\varphi^d$: the battery energy output efficiency.

In reality, electric vehicles consume extra energy when accelerating initially since they have to overcome resistance to change their motion state. In contrast, when electric vehicles brake, they utilize resistance to slow down, and the motor reverses its function to convert some of the mechanical energy into electrical energy, which is then stored in the battery, leading to energy recovery.

Incorporating the electric vehicle start-stop process into the energy consumption model enhances the ability of airport management departments to refine energy consumption predictions and develop more practical electric vehicle dispatching plans. Based on Basso et al. (2019), this paper derives the instantaneous mechanical energy, denoted as ${q_m}^{\uparrow}$ (for acceleration), ${q_m}^{\rightarrow}$ (for constant speed), ${q_m}^{\downarrow}$ (for deceleration), as shown in equations (2) to (4). For a detailed derivation process, please refer to Section 4.4.

$$\boldsymbol{{q_m}^{\uparrow} = (C_r \cdot m \cdot g + \frac{1}{2} \cdot \rho \cdot A \cdot C_d \cdot \frac{v^2}{2} + m \cdot a) \cdot \frac{v^2}{2a} \cdot \emptyset^d \cdot \varphi^d \cdot \epsilon} \tag{2}$$

$$q_m{}^{\rightarrow} = (C_r \cdot m \cdot g + \frac{1}{2} \cdot \rho \cdot A \cdot C_d \cdot v^2) \cdot \left(l - \frac{v^2}{a}\right) \cdot \emptyset^d \cdot \varphi^d \cdot \epsilon \tag{3}$$

$$q_m{}^{\downarrow} = (C_r \cdot m \cdot g + \frac{1}{2} \cdot \rho \cdot A \cdot C_d \cdot \frac{v^2}{2} - m \cdot a) \cdot \frac{v^2}{2a} \cdot \emptyset^r \cdot \varphi^r \cdot \epsilon \tag{4}$$

Where:

$v$: the maximum attainable driving speed of the vehicle (m/s).

$a$: the magnitude of the acceleration during the acceleration and deceleration processes (m/s$^2$).

$\emptyset^r$: the input efficiency of the drive motor when recuperating electrical energy.

$\varphi^r$: the battery energy recuperating efficiency.

Furthermore, we calculate the cumulative energy consumption of the electric vehicle on a hundred-meter road segment, denoted as $q_{\text{cum}}$:

$$q_{cum} = \begin{cases} \left(C_r \cdot m \cdot g + \frac{1}{2} \cdot \rho \cdot A \cdot C_d \cdot \frac{v^2}{2} + m \cdot a\right) \cdot l \cdot \emptyset^d \cdot \varphi^d \cdot \epsilon, 0 \leq l < \frac{v^2}{2a} \\ q_m{}^{\uparrow} + (C_r \cdot m \cdot g + \frac{1}{2} \cdot \rho \cdot A \cdot C_d \cdot v^2) \cdot \left(l - \frac{v^2}{2a}\right) \cdot \emptyset^d \cdot \varphi^d \cdot \epsilon, \frac{v^2}{2a} \leq l < 100 - \frac{v^2}{2a} \\ q_m{}^{\uparrow} + q_m{}^{\rightarrow} + (C_r \cdot m \cdot g + \frac{1}{2} \cdot \rho \cdot A \cdot C_d \cdot \frac{v^2}{2} - m \cdot a) \cdot \left(l - 100 + \frac{v^2}{2a}\right) \cdot \emptyset^r \cdot \varphi^r \cdot \epsilon \\ , 100 - \frac{v^2}{2a} \leq l \leq 100 \end{cases} \tag{5}$$

To provide a more intuitive comparison of the energy consumption results obtained using the traditional model (Equation (1)) and the improved model proposed in this paper (Equation (5)), we compute cumulative energy consumption at various maximum speeds for airport electric tractors, including low (3 m/s), medium (5 m/s), and high (7 m/s) speeds. The results are shown in Figure 1.

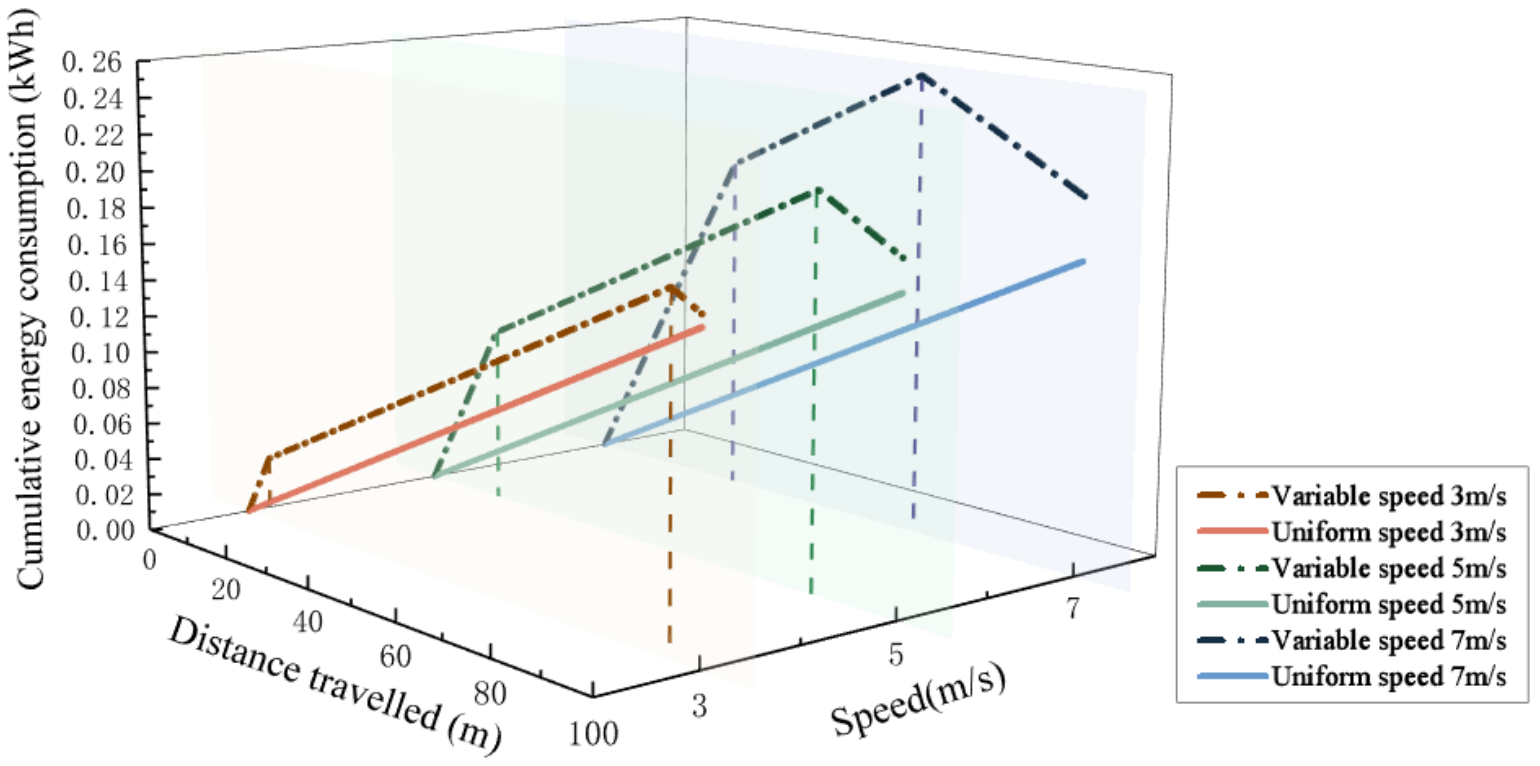

**Fig. 1.** Energy consumption-distance curve within 100 meters.

Fig. 1 shows that under the same speed settings, considering the start-stop process of electric vehicles results in higher energy consumption compared to not considering this process. Furthermore, as the maximum driving speed increases, this difference becomes more significant. When comparing the curve that accounts for the start-stop process to the one that disregards it, it is obvious that the former exhibits a sharper incline during acceleration, a parallel slope during constant speed, and even a gentler descent during deceleration. Energy consumption over a 100-meter distance reveals a noticeable difference when taking the start-stop process into account versus omitting it. This leads to the conclusion that in airport field dispatching problems with many nodes and tight spacing, neglecting the start-stop process can have a more significant impact on the energy consumption estimation.

### *3.2. Piecewise linear charging function*

While the trend of the charging curve is known, expressing it with mathematical equations is complex due to its dependence on factors such as current, voltage, self-recovery, and temperature, as depicted in Fig. 2. The charging function can be divided into two stages based on the battery charge level: the first stage is the constant current charging phase, where the charging current remains constant, and the battery charge level linearly increases over time until it reaches a threshold, at which point the battery charge level is approximately 80% of its capacity. The second stage is the constant voltage charging phase, where, to prevent battery damage, the current decreases while the voltage remains constant. The battery charge level increases synchronously with time during this phase but not in a linear fashion.

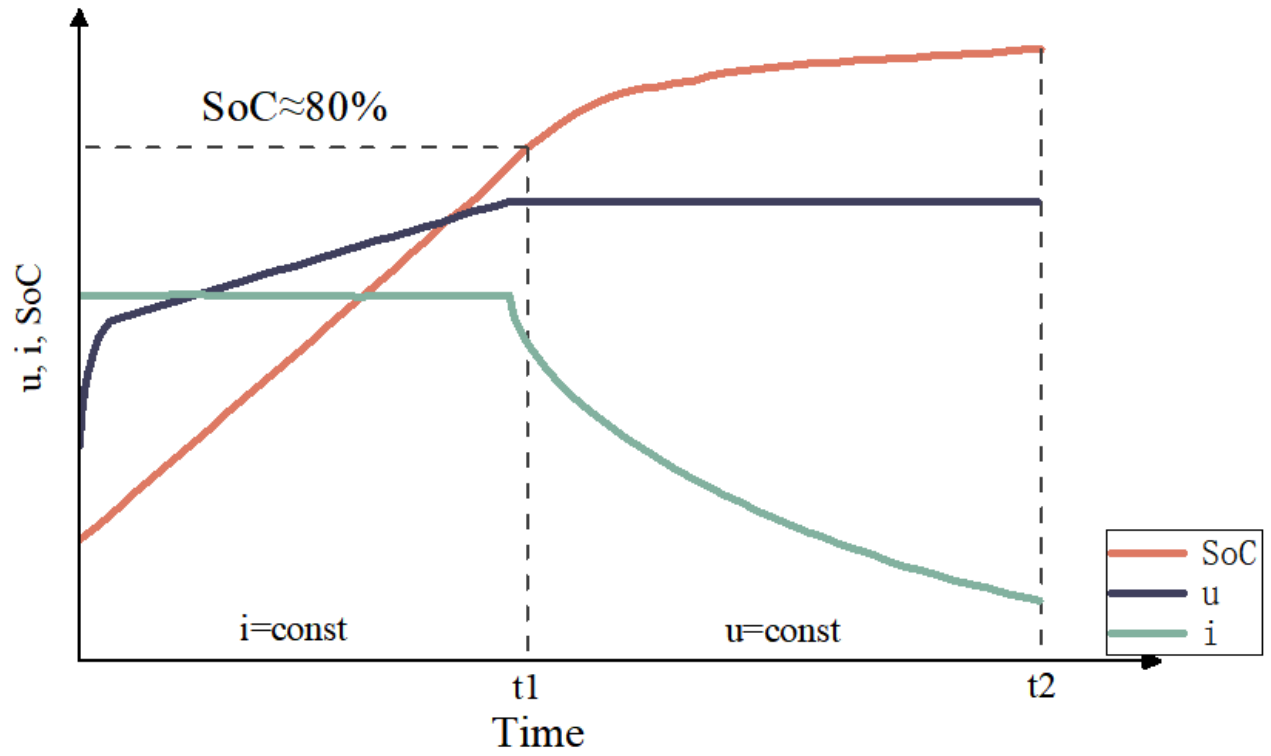

**Fig. 2.** Typical charging curve. (i, u, and SoC represent current, terminal voltage, and State of Charge respectively. SoC is equivalent to the battery charge level.) (Source: Hõimoja et al., 2012)

Previous research often took simplistic fitting approaches, using either optimistic or pessimistic criteria to approximate the charging curve. The optimistic criteria align the charging rate for the entire charging process with the slope of the first phase of the charging function, as shown in Fig. 3(a). Nevertheless, this fitting approach typically underestimates the time required to achieve a full charge. Conversely, the pessimistic criteria set the charging rate for the entire charging process to align with the slope of the line connecting the initial and final observed battery charge levels, as illustrated in Fig. 3(b). However, this fitting approach often results in charging rates that are slower than those encountered in most real-world scenarios.

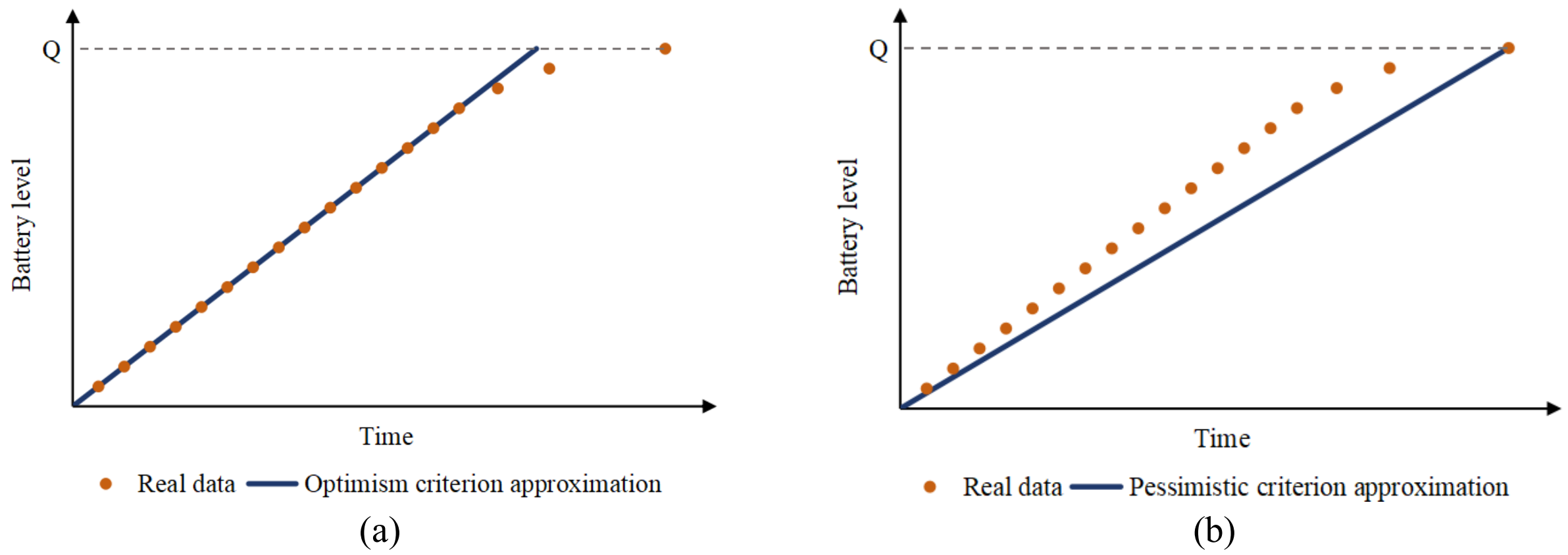

**Fig. 3.** Charging function a) under the optimistic criteria; b) under the pessimistic criteria. (Source: Uhrig et al., 2015)

The two criteria mentioned above are not well-suited to reality. Therefore, this paper adopts a piecewise linear approximation method, as proposed by Montoya et al. (2017). Using real airport charging station data, an equally proportioned non-linear fitting is performed on the charging process to obtain an improved charging function image that closely aligns with reality, as depicted in Fig. 4. The function expression is shown in equation (6).

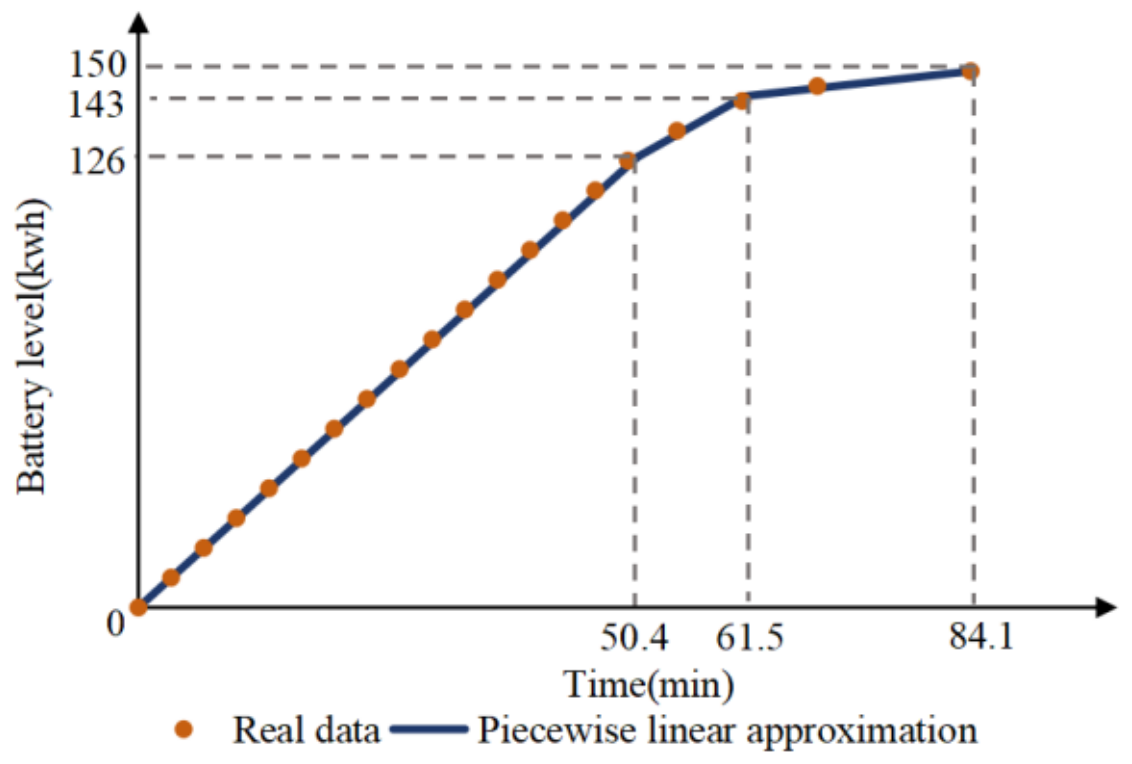

**Fig. 4.** Improved charging function.

$$SoC = \begin{cases} 2.5t, t\epsilon(0, 50.4] \\ 1.5t + 50.4, t\epsilon(50.4,61.5] \\ 0.325t + 122.7, t\epsilon(61.5,84.1] \end{cases} \tag{6}$$

Where:

$SoC$: the battery charge level (kwh).

$t$: the charging time (min).

Utilizing the piecewise linear charging function can help the airport precisely estimate the battery charge level of electric service vehicles, facilitating the development of practical dispatching schedules for these vehicles.

### *3.3. Charge-discharge coupling strategy*

The complete service process of an electric vehicle is divided into two parts: discharging and charging. During the discharging process, the pacing strategy is reflected in the setting of the maximum driving speed, which further determines the energy consumption during the start-stop processes of the electric vehicle. During the charging process, the charging strategy is implemented by configuring the maximum battery charge threshold, which in turn determines the charging duration.

Efficient pacing strategies play a vital role in the operation of airport ground service vehicles. However, there is currently a lack of research dedicated to pacing strategies for these electric vehicles. Drivers frequently depend on personal preferences and experience in their driving practices, which often lack a scientific foundation. Such driving behavior can lead to inefficiencies in operations, energy wastage, and disordered airport management. Meanwhile, charging strategies exhibit diversity in other domains, while airports commonly adopt a full charging approach. The introduction of piecewise linear charging functions has made the optimal charging strategy no longer universally static, especially for airports with different throughputs. Instead, it is essential to adaptively choose the optimal battery charge level (threshold) based on specific circumstances.

Optimal vehicle fleet management is attainable through the integration of pacing and charging strategies, achieved by adjusting maximum speeds and the maximum battery charge threshold. To achieve that, this paper introduces a novel concept of charge-discharge coupling, considering both pacing and charging strategies. This approach aims to optimize the discharging and charging processes in airport electric vehicle operations. According to the Safety Management Rules issued by the Civil Aviation Administration of China as published in CCAR-331SB-R1, ground service vehicles should not exceed a speed of 25 kilometers per hour. Therefore, this study evaluates speeds across five categories (5km/h, 10km/h, 15km/h, 20km/h, and 25km/h) to determine the most effective pacing strategies. Also, we categorize three charging strategies based on the charging level: high level (charging up to 100% of the battery capacity), medium level (charging up to 90% of the battery capacity), and low level (charging up to 80% of the battery capacity). The goal is to design charging strategies that are most suitable. Fig. 5 illustrates a schematic diagram of the charge-discharge coupling strategy.

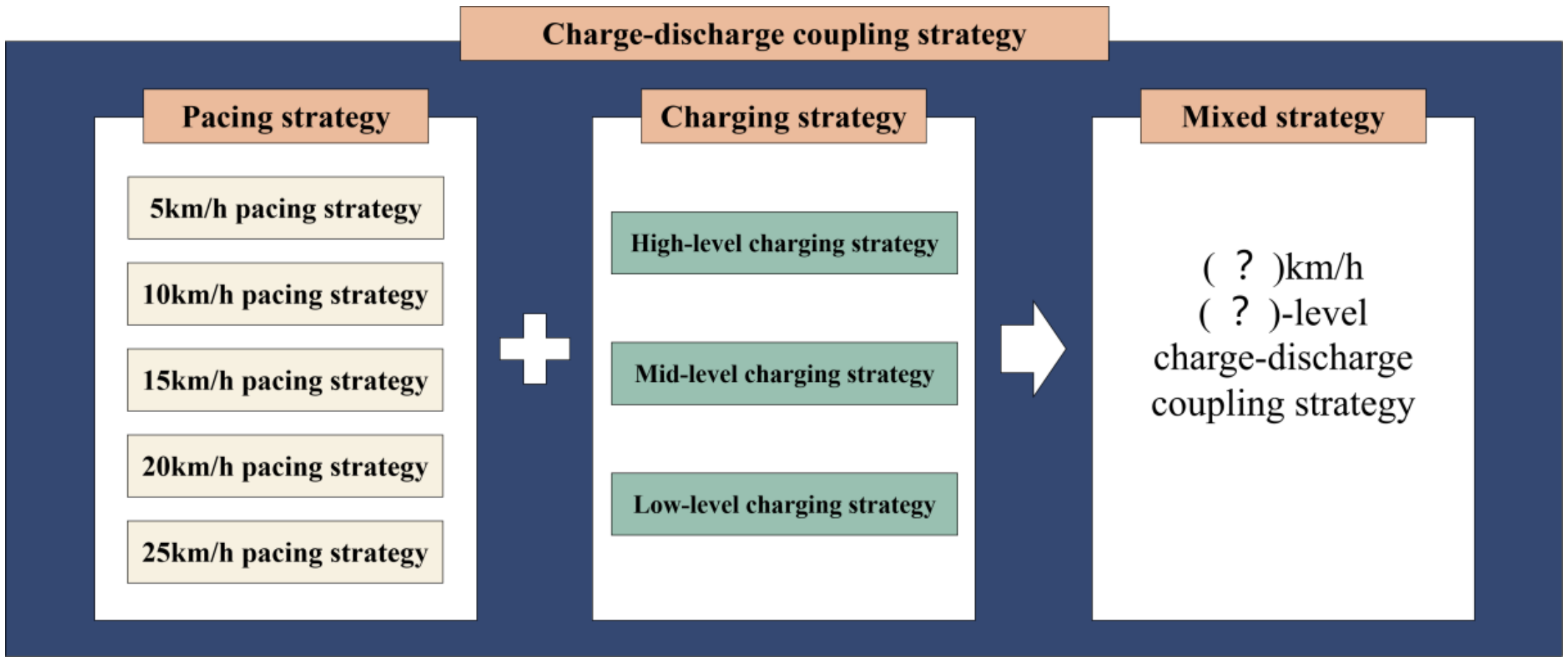

**Fig. 5.** The schematic diagram of the charge-discharge coupling strategy.

## 4. Mathematical model

This section delineates the critical parameters and variables necessary to tackle the dispatching problem. Building upon this foundation, it proposes a dispatching model designed specifically for airport electric tractors. This model incorporates time windows and encompasses a complex energy consumption model.

### *4.1. Parameters*

In the case of flights where the required towing services are predetermined, we operate under the assumption that the airport employs the same type of electric tractors for providing these services. Consider a directed weighted network represented as $G=(S,V)$, where $S$ represents the set of nodes, and $V = \{(i,j)|i,j \in S, i \neq j\}$ is the set of all feasible paths between any pair of nodes. The node set $S$ is defined as $\{d\} \cup P \cup H$, with $d$ denoting the parking lot node where all tractors depart from and eventually return to after completing all service tasks. $P$ signifies the set of aircraft stands. Note that aircraft stands are utilized by multiple flights, and for the sake of mathematical modelling, an aircraft stand serving different flights is labeled as different aircraft stand nodes. $H$ represents the set of charging stations, which can be visited multiple times. Lastly, let $N$ represent the set of electric tractors.

For each aircraft stand $i$, a tractor performs a sequence of towing service and power supply service for aircraft. This paper assumes that the towing distance required for each aircraft is uniform and equal to $w$. When towing the aircraft, the tractor starts from rest and accelerates to a speed $u$ with an acceleration $\varepsilon$. Subsequently, during the towing process, it maintains a constant speed $u$. As it approaches the target parking position (the towing destination), the tractor decelerates with a deceleration $\varepsilon$ and comes to a stop alongside the aircraft. Subsequently, power supply services are initiated. The tractor supplies power to the aircraft for the air conditioning system, lighting system, and engine. The specific service duration and energy consumption per unit of time are contingent upon the aircraft type. Accordingly, the total service time for an aircraft parking at aircraft stand $i$ is denoted as $t_i$ and comprises two components: the towing time and the power supply time. Given that this paper assumes constant values for towing distance, maximum speed, and acceleration, the towing time is uniquely determined. Therefore, the total service time $t_i$ depends solely on the power supply time, which varies based on the aircraft type. The total energy consumption required for aircraft parking at aircraft stand $i$ is represented as $q_i$, and it also consists of two distinct components: energy consumption during towing and power supply, both of which vary depending on the aircraft type.

The service time for each aircraft at aircraft stand $i \in P$ is subject to specific time window constraints. The time that tractor $k$ arrives at aircraft stand $i$ must fall between the earliest arrival time $e_i$ and the latest arrival time $\tau_i$. For each path $(i, j)$, it has path length $l_{ij}$, tractor travel time $t_{ij}$. The energy consumption associated with any tractors travelling along the path $(i, j)$ is denoted by $q_{ij}$. The tractor travel process can be described as follows: the tractor accelerates uniformly with an acceleration $a$ from the origin until it reaches the maximum travel speed $v$ (we will discuss the choice of value of $v$ in the later sections). Subsequently, on the designated path, it maintains a constant speed $v$. When approaching the destination node, the tractor uniformly decelerates with a deceleration $a$ and comes to a stop at the destination node.

Each tractor has a maximum battery capacity of $B$ and a minimum battery threshold of $0.2B$ (airports establish a minimum battery threshold to enhance the rechargeable lifespan and optimize battery performance). During the process of providing towing services for flights, when a tractor $k$ arrives at aircraft stand $i$ with a remaining battery charge level below the minimum threshold, it needs to visit nearby charging stations to charge its battery. The charging time at charging station $i$ is $\delta_i$, where the values of $\delta_i$ are defined in discrete segments according to the actual charging curve. The battery in the tractor can be charged up to the maximum threshold of $\gamma B$, where $\gamma \in (0,1)$ ($\gamma$ represents a coefficient defining the maximum capacity threshold, which imposes constraints on the fixed maximum charging capacity; we will discuss the value of $\gamma$ in the later sections).

In this paper, the following decision variables are defined: if tractor $k$ travels on path $(i, j)$, $x_{ijk}$ is assigned the value of 1; otherwise, it is 0. $s_{ik}$ represents the arrival time of tractor $k$ at aircraft stand $i$. $p_{ik}$ represents the remaining battery charge level when tractor $k$ arrives at node $i \in S$. Table 1 summarizes the parameters and variables.

**Table 1**

Parameter and variable definitions.

| Sets | |
|---|---|
| $\{d\}$ | Parking lot node |
| $P$ | Airport stand nodes set. The same stand where different flights are parked is considered to be different stand nodes |
| $H$ | Charging station nodes set |
| $S$ | Nodes set, $S = \{d\} \cup P \cup H$ |
| $V$ | Routes set, $V = \{(i,j) \mid i,j \in S, i \neq j\}$ |
| $N$ | Tractors set, $k \in N$ |
| $G$ | Airport directed weighed network, $G=(S,V)$ |
| Parameters | |
| $l_{ij}$ | Distance between node $i$ and $j$, $i,j \in S$ (m) |
| $w$ | Distance of towing service (m) |
| $t_{ij}$ | Travel time between node $i$ and $j$, $i,j \in S$ (min) |
| $t_i$ | Service time at stand node $i$, consisting of towing time and power supply time (min) |
| $v$ | Maximum allowable travel speed (m/s) |
| $u$ | Maximum allowable towing speed (m/s) |
| $a$ | Magnitude of acceleration during the travel process (m/s$^2$) |
| $\varepsilon$ | Magnitude of acceleration during the towing process (m/s$^2$) |
| $\delta_i$ | Charging time at charging stations $i$. The values of $\delta_i$ are defined in discrete segments according to the actual charging curve (min) |
| $B$ | Battery capacity (kwh) |
| $\gamma$ | Maximum battery capacity threshold coefficient |
| $q_i$ | Service energy consumption at stand node $i$, consisting of towing energy consumption and power supply energy consumption (kwh) |
| $q_{ij}$ | Travelling energy consumption between node $i$ and $j$, $i,j \in S$ (kwh) |
| $e_i$ | Earliest service time required by the aircraft parked at the stand node $i$ (min) |
| $\tau_i$ | Latest service time required by the aircraft parked at the stand node $i$ (min) |
| Variables | |
| $x_{ijk}$ | Decision variable, whether $(i,j)$ is passed by tractor $k$ |
| $s_{ik}$ | Arrival time of tractor $k$ at stand node $i$, if tractor $k$ serves the flight located at this node |
| $p_{ik}$ | Remaining power of tractor $k$ arriving at node $i$, $i \in S$, if tractor $k$ reaches this node |

### *4.2. Assumptions*

The primary focus of the model is the optimization of trajectories for airport electric tractors under various charge-discharge coupling strategies. This model includes tractor start-stop procedures, integrates flight services, and accommodates a piecewise linear charging function, which contributes to the complexity of the constraints. To streamline the model, the following assumptions are made in advance:

1. Each flight only needs the assistance of a single tractor for towing, and once the service starts, it proceeds without interruption until completion.
2. The model does not consider the influence of traffic conditions. Once the tractor reaches its maximum speed, it maintains that speed without regard to traffic conditions.
3. Charging stations have an ample supply of charging piles, allowing multiple tractors to charge simultaneously without queuing.
4. All electric tractors share the same model and performance characteristics. Tractors start their journeys with a full battery charge and return to the parking lot after completing all service tasks.

### *4.3. Simple time network model*

Flight ground service procedures refer to a series of ground services provided to an aircraft after it lands at the aircraft stand. The most crucial procedures include guidance, baggage handling, refueling, cleaning, catering, and so on. There are sequential constraints among these servicing tasks. For instance, for transit flights, the guidance service should occur before all other services. Each service also has its own time constraints. For instance, for flights with transit times of less than 80 minutes, the duration of the guidance service should fall within a time window of 5 to 10 minutes.

Considering that originating flights typically receive their ground services on the day of arrival (the day before their departure from the airport) and have longer transit times with ample service windows, this paper will focus exclusively on the scenarios involving transit flights. As the initial step in the servicing procedures for transit flights, the dispatching of towing services directly impacts the smooth execution of subsequent services. Therefore, it is essential to comprehensively consider the sequence and timing constraints of all flight service procedures to ensure that the assigned towing service time windows can guide actual operations. The advantage of the Simple Time Network (STN) model lies in its efficiency in handling various types of time window constraints, including those that are elastic and uncertain. Thus, this study employs the STN model to determine the towing service time windows. The specific steps for solving this are as follows:

Step 1: Clearly define the flight ground service procedures and analyze the types of ground service vehicles required for transit flights.

Step 2: To represent the various time constraints for different service tasks based on the varying transit times for flights, we can use an interval notation formalized as $Q = (X, C)$. Here, $X$ is the set of time variables, denoted as $X = \{x_0, x_1, \ldots, x_n\}$, where each element in set $X$ represents the start or end time of various services. The set $C$ comprises the relationships between these time variables, and it includes constraints of the form $c_{ij}\colon x_j - x_i \leq b_{ij}$, where $c_{ij}$ signifies that $x_j$ occurs after $x_i$ by a duration of $b_{ij}$ time units. Here, $b_{ij} \in Z$, where $Z$ represents the set of integers.

Step 3: Integrate the sequential constraints with the timing constraints to construct a STN Model for flight servicing procedures. The STN is a directed constrained network graph, denoted as $R = (J, L)$. Here, the vertex set $J$ comprises all the points corresponding to time variables, and $L$ represents the set of edges between any pair of two points in the set $J$. If the constraints $c_{ij}: x_j - x_i \leq b_{ij}$ are covered in $C$, then add an edge from $i$ to $j$ in $L$, assigning it a weight of $b_{ij}$.

Step 4: Perform decoupling operations on the STN to obtain the earliest start time, earliest end time, latest start time, and latest end time for the towing service of each flight.

### *4.4. Energy consumption model*

The energy consumption of electric vehicles entails the conversion of electrical energy stored in the battery into electrical energy for the propulsion motor, subsequently transformed into mechanical energy to propel the vehicle forward. Given that the consumption of mechanical energy is easier to measure and assess (e.g., vehicle resistance), we choose to analyze energy consumption in the reverse manner. Therefore, the first step is to calculate the instantaneous mechanical power $P_m(t)$ of the tractor at time $t$:

$$F_m(t) = F_r(t) + F_a(t) + F_g(t) + m \cdot a(t) \tag{7}$$

Where:

$F_m(t)$: the force on the powertrain system at time $t$ (N).

$F_r(t)$: the rotational resistance acting on the tractor at time $t$ (N).

$F_a(t)$: the air resistance acting on the tractor at time $t$ (N).

$F_g(t)$: the component of gravity acting on the tractor in the direction parallel to the road surface at time $t$ (N).

$m$: the mass of the tractor (kg).

$a(t)$: the acceleration of the tractor at time $t$ (m/s$^2$).

Here are the formulas for $F_r(t)$, $F_a(t)$, and $F_g(t)$:

$$\boldsymbol{F_r(t) = C_r \cdot m \cdot g \cdot cos\,\alpha(t)} \tag{8}$$

$$\boldsymbol{F_a(t) = \frac{1}{2} \cdot \rho \cdot A \cdot C_d \cdot v(t)^2} \tag{9}$$

$$\boldsymbol{F_g(t) = m \cdot g \cdot \sin \alpha(t)} \tag{10}$$

Where:

$C_r$: the coefficient of rolling friction.

$g$: the acceleration due to gravity (N/kg).

$\alpha(t)$: the road gradient at time $t$ (m/s$^2$).

$\rho$: air density (kg/m$^3$).

$A$: the frontal area of the tractor (m$^2$).

$C_d$: the coefficient of air resistance.

$v(t)$: the instantaneous speed of the tractor at time $t$ (m/s).

According to the Civil Airport Flight Area Technical Standards, road gradients are specified to be below 1%. Therefore, $\sin\alpha(t) \approx 0, \cos\alpha(t) \approx 1$. The formula for calculating the instantaneous mechanical power $P_m(t)$ is as shown in Equation (11):

$$\boldsymbol{P_m(t)=C_r \cdot m \cdot g \cdot v(t) + \frac{1}{2} \cdot \rho \cdot A \cdot C_d \cdot v(t)^3 + m \cdot a(t) \cdot v(t)} \tag{11}$$

The second step involves reversing the calculation process to determine the electric power, denoted as $P_e(t)$, delivered by the drive motor at time $t$. This paper uses the conversion formula between $P_m(t)$ and $P_e(t)$ as proposed by Goeke (2015). They directly converted torque and rotational speed into mechanical power and established the relationship between $P_m(t)$ and the electric energy output (or recuperating energy) $P_e(t)$ by conducting a linear regression analysis with a y-axis intercept of zero on the converted values:

$$\boldsymbol{P_e(t)=\begin{cases} \emptyset^d \cdot P_m(t), P_m(t) \geq 0 \\ \emptyset^r \cdot P_m(t), P_m(t) < 0 \end{cases}} \tag{12}$$

Where:

$\emptyset^d$: the output efficiency of the drive motor when releasing electrical energy.

$\emptyset^r$: the input efficiency of the drive motor when recuperating electrical energy.

The third step involves reversing the calculation process to determine the energy $P_q(t)$ that the battery needs to release at time $t$. This paper also utilizes the conversion formula between $P_e(t)$ and $P_q(t)$ as proposed by Goeke (2015). They considered external factors such as environmental temperature and the current battery charge level, and used linear regression to establish the relationship between $P_e(t)$ and battery discharge (or recuperate) $P_q(t)$:

$$\boldsymbol{P_q(t)=\begin{cases} \varphi^d \cdot P_e(t), P_e(t) \geq 0 \\ \varphi^r \cdot P_e(t), P_e(t) < 0 \end{cases}} \tag{13}$$

Where:

$\varphi^d$: the battery energy output efficiency.

$\varphi^r$: the battery energy recuperating efficiency.

In the fourth step, calculate the energy consumption $q_{ij}$ of the tractor travelling on path $(i,j)$. This paper refers to the energy consumption integral formula proposed by Basso et al.(2019). The towing process on path $(i,j)$ is divided into three stages: accelerating, maintaining maximum speed, and decelerating, with energy consumption represented by $q_{ij}^{\uparrow}$, $q_{ij}^{\rightarrow}$, and $q_{ij}^{\downarrow}$ respectively. Assuming an initial speed of 0 for the tractor, a maximum allowable travel speed of $v$, and acceleration and deceleration magnitudes of $a$ travelling a distance of $l_{ij}$, the following equations apply:

$$\boldsymbol{q_{ij}}^{\uparrow} = (\boldsymbol{C_r} \cdot \boldsymbol{m} \cdot \boldsymbol{g} + \frac{1}{2} \cdot \boldsymbol{\rho} \cdot \boldsymbol{A} \cdot \boldsymbol{C_d} \cdot \frac{\boldsymbol{v}^2}{2} + \boldsymbol{m} \cdot \boldsymbol{a}) \cdot \frac{\boldsymbol{v}^2}{2\boldsymbol{a}} \cdot \emptyset^{\boldsymbol{d}} \cdot \boldsymbol{\varphi}^{\boldsymbol{d}} \cdot \epsilon \tag{14}$$

$$\boldsymbol{q_{ij}}^{\rightarrow} = (\boldsymbol{C_r} \cdot \boldsymbol{m} \cdot \boldsymbol{g} + \frac{1}{2} \cdot \boldsymbol{\rho} \cdot \boldsymbol{A} \cdot \boldsymbol{C_d} \cdot \boldsymbol{v}^2) \cdot \left(\boldsymbol{l_{ij}} - \frac{\boldsymbol{v}^2}{\boldsymbol{a}}\right) \cdot \emptyset^{\boldsymbol{d}} \cdot \boldsymbol{\varphi}^{\boldsymbol{d}} \cdot \epsilon \tag{15}$$

$$\boldsymbol{q_{ij}}^{\downarrow} = (\boldsymbol{C_r} \cdot \boldsymbol{m} \cdot \boldsymbol{g} + \frac{1}{2} \cdot \boldsymbol{\rho} \cdot \boldsymbol{A} \cdot \boldsymbol{C_d} \cdot \frac{\boldsymbol{v}^2}{2} - \boldsymbol{m} \cdot \boldsymbol{a}) \cdot \frac{\boldsymbol{v}^2}{2\boldsymbol{a}} \cdot \emptyset^{\boldsymbol{r}} \cdot \boldsymbol{\varphi}^{\boldsymbol{r}} \cdot \epsilon \tag{16}$$

$$\boldsymbol{q_{ij}} = \boldsymbol{q_{ij}}^{\uparrow} + \boldsymbol{q_{ij}}^{\rightarrow} + \boldsymbol{q_{ij}}^{\downarrow} \tag{17}$$

The service energy consumption of electric tractors can be subdivided into towing energy consumption and auxiliary power unit (APU) replacement energy consumption (van Baaren, 2019). Thus, the service energy consumption of the tractor at aircraft stand $i$ is denoted as $q_i$, and it consists of two parts: towing energy consumption $q_i^{DRAG}$ and APU replacement energy consumption $q_i^{APU}$. Specifically, the calculation of towing energy consumption $q_i^{DRAG}$ is the same as it for travelling energy consumption. The power consumption during acceleration, constant speed, and deceleration phases is represented by $q_i^{DRAG^{\uparrow}}$, $q_i^{DRAG^{\rightarrow}}$, and $q_i^{DRAG^{\downarrow}}$ respectively. The difference is that when the tractor provides aircraft towing services, we need to consider the additional mass of the towed aircraft $M_i$. After the towing service is completed, the tractor provides power to the aircraft for air conditioning system, lighting system, and engine start, with corresponding energy consumption rates of $q_i^{AIR}$, $q_i^{ILL}$, and $q_i^{LAU}$, and the required service times are $t_i^{AIR}$, $t_i^{ILL}$, and $t_i^{LAU}$. For each aircraft stand $i$, we have the following formulas:

$$\boldsymbol{q_i^{DRAG^{\uparrow}}} = \left[\boldsymbol{C_r} \cdot (\boldsymbol{m} + \boldsymbol{M_i}) \cdot \boldsymbol{g} + \frac{1}{2} \cdot \boldsymbol{\rho} \cdot \boldsymbol{A} \cdot \boldsymbol{C_d} \cdot \frac{\boldsymbol{u}^2}{2} + (\boldsymbol{m} + \boldsymbol{M_i}) \cdot \boldsymbol{\varepsilon}\right] \cdot \frac{\boldsymbol{u}^2}{2\boldsymbol{\varepsilon}} \cdot \emptyset^{\boldsymbol{d}} \cdot \boldsymbol{\varphi}^{\boldsymbol{d}} \cdot \epsilon \tag{18}$$

$$\boldsymbol{q_i^{DRAG^{\rightarrow}}} = \left[\boldsymbol{C_r} \cdot (\boldsymbol{m} + \boldsymbol{M_i}) \cdot \boldsymbol{g} + \frac{1}{2} \cdot \boldsymbol{\rho} \cdot \boldsymbol{A} \cdot \boldsymbol{C_d} \cdot \boldsymbol{u}^2\right] \cdot \left(\boldsymbol{w} - \frac{\boldsymbol{u}^2}{\boldsymbol{\varepsilon}}\right) \cdot \emptyset^{\boldsymbol{d}} \cdot \boldsymbol{\varphi}^{\boldsymbol{d}} \cdot \epsilon \tag{19}$$

$$\boldsymbol{q_i^{DRAG^{\downarrow}}} = \left[\boldsymbol{C_r} \cdot (\boldsymbol{m} + \boldsymbol{M_i}) \cdot \boldsymbol{g} + \frac{1}{2} \cdot \boldsymbol{\rho} \cdot \boldsymbol{A} \cdot \boldsymbol{C_d} \cdot \frac{\boldsymbol{u}^2}{2} - (\boldsymbol{m} + \boldsymbol{M_i}) \cdot \boldsymbol{\varepsilon}\right] \cdot \frac{\boldsymbol{u}^2}{2\boldsymbol{\varepsilon}} \cdot \emptyset^{\boldsymbol{r}} \cdot \boldsymbol{\varphi}^{\boldsymbol{r}} \cdot \epsilon \tag{20}$$

$$\boldsymbol{q_i^{DRAG}} = \boldsymbol{q_i^{DRAG^{\uparrow}}} + \boldsymbol{q_i^{DRAG^{\rightarrow}}} + \boldsymbol{q_i^{DRAG^{\downarrow}}} \tag{21}$$

$$\boldsymbol{q_i^{APU} = q_i^{AIR} * t_i^{AIR} + q_i^{ILL} * t_i^{ILL} + q_i^{LAU} * t_i^{LAU}} \tag{22}$$

$$\boldsymbol{q_i = q_i^{DRAG} + q_i^{APU}} \tag{23}$$

*4.5. Dispatching model*

This paper, from the perspective of airport operational efficiency, establishes an objective function for the Electric Vehicle Routing Problem with Time Windows (EVRPTW-PR) model, which aims to minimize the total cost. The objective function primarily takes into account fixed costs, charging costs, maintenance costs, and time window penalty costs.

The fixed cost component $f_1$ primarily includes tractor depreciation expenses, driver labor costs, meal expenses, and so on. In this paper, it is assumed that the fixed cost per tractor is the same, so $f_1$ is only dependent on the total number of tractors. The calculation formula is as follows:

$$f_1{=}c_1 \sum_{j\epsilon P} \sum_{k\in N} x_{djk} \tag{24}$$

Where $c_1$ represents the unit tractor fixed cost.

The charging cost $f_2$ is related to the energy consumption. Each time a tractor leaves the charging station, the charging cost is calculated based on a single charge amount of $\gamma B - p_{ik}$. Each time a tractor returns to the parking lot, the charging cost is calculated based on the used electricity of $B - p_{dk}$. The calculation formula is as follows:

$$f_2{=}c_2\left[\sum_{i\epsilon H} \sum_{j\epsilon P} \sum_{k\in N} (\gamma B - p_{ik}) \cdot x_{ijk} + \sum_{j\epsilon P} \sum_{k\in N} (B - p_{dk}) \cdot x_{jdk}\right] \tag{25}$$

Where $c_2$ represents the unit charging cost.

The maintenance cost $f_3$ mainly includes the maintenance fees for both the tractor and the battery. This paper assumes that these two costs are linearly related to the travelling mileage. The calculation formula is as follows:

$$f_3{=}c_3 \sum_{i\epsilon S} \sum_{j\in S, i!=j} \sum_{k\in N} l_{ij} \cdot x_{ijk} \tag{26}$$

Where $c_3$ represents the unit mileage maintenance cost.

The time window penalty cost $f_4$ is applied when the tractor arrives either before or after the designated service time window. If the tractor arrives ahead of the service time window, it incurs waiting time equal to $e_j - s_{jk}$. Conversely, if the tractor arrives after the service time window, it is considered a delay, with the delay time being $s_{jk} - \tau_j$. The calculation formula is as follows:

$$f_4{=}\sum_{i\epsilon S} \sum_{j\in P, i!=j} \sum_{k\in N} \left[c_e \cdot max(e_j - s_{jk}, 0) + c_\tau \cdot max(s_{jk} - \tau_j, 0)\right] \cdot x_{ijk} \tag{27}$$

Where $c_e$ represents the cost per unit of waiting time, and $c_\tau$ represents the cost per unit of delay time.

In summary, the objective function $f$ that minimizes the total cost in this paper is as follows:

$$Minimize\ f = f_1 + f_2 + f_3 + f_4 \tag{28}$$

The constraint set encompasses sequential constraints, time constraints, and energy constraints. The foundational model is constructed as follows:

$$\sum_{i\epsilon S, i!=j} \sum_{k\in N} x_{ijk} = 1, \forall j \in P \tag{29}$$

$$\sum_{j\epsilon P} x_{djk} = 1, \forall k \in N \tag{30}$$

$$\sum_{j\epsilon S, i!=j} x_{ijk} - \sum_{j\epsilon S, i!=j} x_{jik} = 0, \forall i \in P \cup H, \forall k \in N \tag{31}$$

$$\sum_{j\epsilon P} x_{jdk} = 1, \forall k \in N \tag{32}$$

$$e_j \leq \sum_{i\epsilon S, i!=j} \sum_{k\epsilon N} s_{jk} \cdot x_{ijk} \leq \tau_j, j \in P \tag{33}$$

$$t_i = \frac{1}{60}\left(\frac{w}{u} + \frac{u}{\varepsilon}\right) + t_i^{AIR} + t_i^{ILL} + t_i^{LAU}, \forall i \in P \tag{34}$$

$$t_{ij} = \frac{1}{60}\left(\frac{l_{ij}}{v} + \frac{v}{a}\right), \forall i \in S, \forall j \in S, i! = j \tag{35}$$

$$s_{ik} + \left(t_i + t_{ij}\right)x_{ijk} - \tau_i\left(1 - x_{ijk}\right) \leq s_{jk}, \forall i \in P, \forall j \in S, i! = j, \forall\ k \in N \tag{36}$$

$$0.2B \leq p_{ik} \leq B, \forall i \in S, \forall k \in N \tag{37}$$

$$\delta_i = \begin{cases} \frac{\gamma B - p_{ik}}{2.5}, \gamma\epsilon(0.2, 0.84] \\ \frac{\gamma B - 0.84B}{1.5} + \frac{0.84B - p_{ik}}{2.5}, \gamma\epsilon(0.84, 0.95] \\ \frac{\gamma B - 0.95B}{0.325} + \frac{0.11B}{1.5} + \frac{0.84B - p_{ik}}{2.5}, \gamma\epsilon(0.95, 1] \end{cases} \quad, \forall i \in H, \forall k \in N \tag{38}$$

$$s_{ik} + \left(\delta_i + t_{ij}\right)x_{ijk} - \tau_i\left(1 - x_{ijk}\right) \leq s_{jk}, \forall i \in H, \forall j \in S, i! = j, \forall k \in N \tag{39}$$

$$q_i^{DRAG^{\uparrow}} = \left[C_r \cdot (m + M_i) \cdot g + \frac{1}{2} \cdot \rho \cdot A \cdot C_d \cdot \frac{u^2}{2} + (m + M_i) \cdot \varepsilon\right] \cdot \frac{u^2}{2\varepsilon} \cdot \emptyset^d \cdot \varphi^d \cdot \epsilon, \forall i \in P \tag{40}$$

$$q_i^{DRAG^{\rightarrow}} = \left[C_r \cdot (m + M_i) \cdot g + \frac{1}{2} \cdot \rho \cdot A \cdot C_d \cdot u^2\right] \cdot \left(w - \frac{u^2}{\varepsilon}\right) \cdot \emptyset^d \cdot \varphi^d \cdot \epsilon, \forall i \in P \tag{41}$$

$$q_i^{DRAG^{\downarrow}} = \left[C_r \cdot (m + M_i) \cdot g + \frac{1}{2} \cdot \rho \cdot A \cdot C_d \cdot \frac{u^2}{2} - (m + M_i) \cdot \varepsilon\right] \cdot \frac{u^2}{2\varepsilon} \cdot \emptyset^d \cdot \varphi^r \cdot \epsilon, \forall i \in P \tag{42}$$

$$q_i = q_i^{DRAG^{\uparrow}} + q_i^{DRAG^{\rightarrow}} + q_i^{DRAG^{\downarrow}} + q_i^{AIR} * t_i^{AIR} + q_i^{ILL} * t_i^{ILL} + q_i^{LAU} * t_i^{LAU}, \forall i \in P \tag{43}$$

$$q_{ij}^{\uparrow} = (C_r \cdot m \cdot g + \frac{1}{2} \cdot \rho \cdot A \cdot C_d \cdot \frac{v^2}{2} + m \cdot a) \cdot \frac{v^2}{2a} \cdot \emptyset^d \cdot \varphi^d \cdot \epsilon, \forall i \in S, \forall j \in S, i! = j \tag{44}$$

$$q_{ij}^{\rightarrow} = (C_r \cdot m \cdot g + \frac{1}{2} \cdot \rho \cdot A \cdot C_d \cdot v^2) \cdot \left(l_{ij} - \frac{v^2}{a}\right) \cdot \emptyset^d \cdot \varphi^d \cdot \epsilon, \forall i \in S, \forall j \in S, i! = j \tag{45}$$

$$q_{ij}^{\downarrow} = (C_r \cdot m \cdot g + \frac{1}{2} \cdot \rho \cdot A \cdot C_d \cdot \frac{v^2}{2} - m \cdot a) \cdot \frac{v^2}{2a} \cdot \emptyset^d \cdot \varphi^r \cdot \epsilon, \forall i \in S, \forall j \in S, i! = j \quad (46)$$

$$q_{ij} = q_{ij}^{\uparrow} + q_{ij}^{\rightarrow} + q_{ij}^{\downarrow}, \forall i \in S, \forall j \in S, i! = j \quad (47)$$

$$p_{jk} \leq p_{ik} - (q_i + q_{ij})x_{ijk} + B(1 - x_{ijk}), \forall i \in P, \forall j \in S, i! = j, \forall k \in N \quad (48)$$

$$p_{jk} \leq \gamma B - q_{ij} \cdot x_{ijk} + B(1 - x_{ijk}), \forall i \in H, \forall j \in S, i! = j, \forall k \in N \quad (49)$$

$$x_{ijk} \in \{0,1\}, \forall i, j \in S, i! = j, k \in N \quad (50)$$

Equation (29) represents that each flight accepts towing service from only one tractor for once. Equation (30) signifies that all tractors depart from the parking node to perform tasks. Equation (31) indicates that when a tractor enters a node, it will eventually depart from that same node. Equation (32) portrays that tractors return to the parking node after completing all service tasks. Equation (33) highlights that all flights must receive traction service within their designated time windows. Equation (34) represents the service time of tractors at aircraft stand $i$. Equation (35) denotes the travel time of tractors on path $(i, j)$. Equation (36) illustrates the time relationship from aircraft stand $i$ to the next node $j$, ensuring that the arrival time of the tractor at the next node $j$ is not less than the sum of the arrival time at aircraft stand $i$, service time at aircraft stand $i$, and travel time on path $(i, j)$. Equation (37) asserts that the remaining available battery charge level of any tractor at any moment is not less than the threshold of $0.2B$ and does not exceed the battery capacity B. Equations (38) represent the charging time of the tractor at charging station $i$. Equation (39) expresses the time relationship from charging station $i$ to the next node $j$, ensuring that the arrival time of the tractor at the next node $j$ is not less than the sum of the arrival time at charging station node $i$, charging time at charging station $i$, and travel time on arc $(i, j)$. Equations (40)-(43) represent the service energy consumption of the tractor at aircraft stand $i$. Equations (44)-(47) represent the travelling energy consumption of the tractor on path $(i, j)$. Equation (48) establishes the energy relationship from aircraft stand $i$ to the next node $j$, ensuring that the tractor's battery charge level at the next node $j$ is not greater than the battery charge level at aircraft stand $i$ minus the energy consumed at aircraft stand $i$ and on path $(i, j)$. Equation (49) defines the energy relationship from charging station node $i$ to the next node $j$, ensuring that the tractor's battery charge level at the next node $j$ is not greater than the battery level charged at charging station $i$ minus the energy consumed on arc $(i, j)$. Equation (50) represents the 0-1 decision variable constraints.

## 5. Algorithm

Considering the complicated energy consumption processes and piecewise linear charging function, the airport tractor dispatch problem is a typical NP-hard problem that is difficult to solve with exact algorithms. Genetic algorithm, as one of the most classic heuristic algorithms, possesses strong global search capabilities and good adaptability, making it widely applicable in the field of trajectory planning and tractor dispatching. It starts by generating initial solutions, retains the population by selecting operators and applies crossover or mutation operations to create new populations. This process is repeated to optimize locally until relatively optimal feasible solutions are found in the solution space.

However, the traditional genetic algorithm is not applicable to all problems, especially the problems with high differentiability. And the energy consumption model constructed in this paper considers the start-stop process, which is highly nonlinear and discrete, making it difficult to be solved only by traditional genetic algorithms. Not only that, Macrina (2019) concluded that the traditional genetic algorithm has two drawbacks: easy to fall into the local optimum and slow convergence of the algorithm. In this regard, this paper incorporates the adaptive strategy proposed by Srinivas M in the crossover and mutation links, and makes adaptive adjustments through iterations. In addition, this paper considers the greedy idea in the initialization of the population, using its local optimization characteristics to iteratively optimize the solution through greedy inverse operation and random search.

### *5.1. Overview*

**(1) Chromosome structure design**

In this paper, chromosomes are defined as the driving trajectories of tractors, encoded using natural numbers. We assume that the length of each chromosome in the population is 2n+1, where $n$ represents the number of aircraft stands, and $m$ represents the number of charging stations. Parking lot is assigned the number 0, aircraft stands are assigned numbers 1, 2, 3, ..., $n$, and charging stations are assigned numbers $n+1, n+2, n+3, \ldots, n+m$. Two instances of parking lot number 0 are inserted at the starting and ending points, thus forming a complete chromosome. For example, assume that there is one charging station and five flight points to be serviced on the apron area. If a chromosome is represented as shown in Fig. 6, then the corresponding driving trajectory for the tractor can be described as follows: The first tractor departs from the parking lot, sequentially serves aircraft stands 2 and 4, and finally returns to the parking lot. The second tractor departs from the parking lot, sequentially serves aircraft stands 1 and 3, then enters charging station 6 for recharging, followed by serving aircraft stand 5, and finally returns to the parking lot.

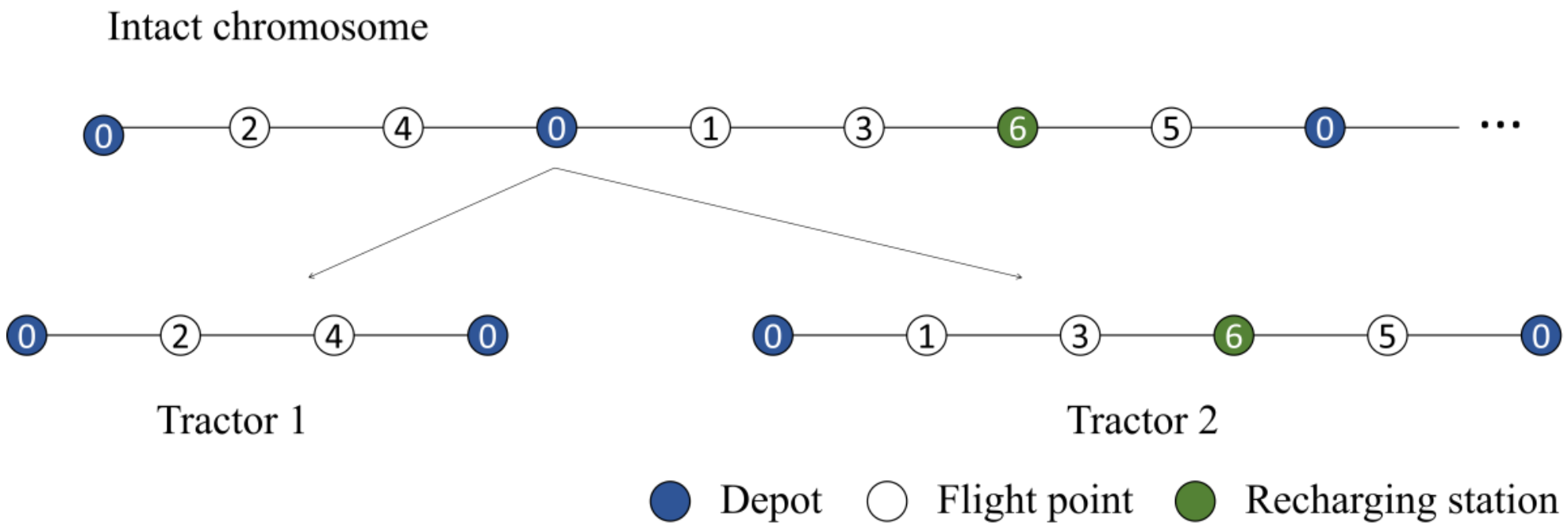

**Fig. 6.** Chromosome decoding process.

**(2) Population initialization**

Typically, population initialization is carried out using a random generation approach. However, this approach often results in individuals with low fitness and a widely dispersed population, which in turn can lead to issues such as a significant difference between feasible solutions and optimal solutions. Therefore, in this paper, the concept of the greedy algorithm is taken into account during the population initialization with the aim of leveraging its local optimization characteristics to enhance the quality of initial solutions. To get a better initial solution, the algorithm used is depicted as follows:

Step 1: Sorting the flights according to the earliest service start times in ascending order.

Step 2: Within the parking lot, assume there is a single tractor $i$ ($i$=1) with a full charge and a time status of 0:00.

Step 3: Choose the unserviced flight $j$ whose earliest service start time is closest to the current time.

Step 4: Check if the battery level for the tractor $i$ meets the constraint. If the battery level is below the minimum threshold of $0.2B$, select the nearest charging station for recharging. Add the charging station to the trajectory of tractor $i$, update its battery level, time status, and return to Step 3. Otherwise, proceed to Step 5.

Step 5: Check if tractor $i$ can catch up with the latest service time of flight $j$. If the tractor $i$ can reach the corresponding aircraft stand within an acceptable delay time, add this stand to the trajectory of tractor $i$, exclude flight $j$ from the set of unserviced flights, update the tractor's battery level and time status, and proceed to Step 6. Otherwise, let tractor $i$ returns to the parking lot. Add another tractor (increment the value of $i$ by 1) to the parking lot, initialize its battery level, time status, and return to Step 3.

Step 6: Check if there are any unserviced flights left. If not, continue to Step 7. Otherwise, return to Step 3.

Step 7: Let tractor $i$ return to the parking lot. Output the assigned service trajectory for all tractors. Calculate the total cost of the trajectory.

**(3) Fitness evaluation**

In this paper, the objective function is defined as minimizing the sum of fixed costs, charging costs, maintenance costs, and time window penalty costs. Since higher fitness values are preferred, in the airport tractor dispatch problem, there is a negative correlation between the objective function value and fitness value. The reciprocal of the objective function is used as the fitness function:

$$Fitness(x) = \frac{1}{Minimize\ f} \tag{51}$$

**(4) Population update**

**a) Selection**

In this study, a combination of the elite strategy and roulette wheel selection is used for population selection. The elite strategy, also referred to as the best preservation strategy, entails transferring the chromosome with the lowest total cost generated in each iteration to the next iteration, thereby accelerating the convergence speed of the genetic algorithm. Meanwhile, the roulette wheel selection introduces the concept of "cumulative probability," influencing the random selection of aircraft stands with higher transition probabilities for visitation, thus increasing the randomness in the trajectory search process.

**b) Crossover**

In this study, an adaptive crossover strategy is used to adjust the crossover probability, denoted as $P_c$, based on the fitness values, $f_a$and $f_b$, of two individuals, $a$ and $b$, participating in the crossover operation. Here, $f$ represents the maximum fitness value among the two individuals involved in the crossover operation, i.e., $f = \max(f_a, f_b)$. When the population tends towards a local optimum, the value of $P_c$ will increase. When the population is more dispersed, $P_c$ decreases. The expression for adjusting the adaptive crossover probability is given as shown in Equation (52):

$$P_c = \begin{cases} P_{c_{max}} - \frac{\left(P_{c_{max}} - P_{c_{min}}\right)(f - f_{avg})}{f_{max} - f_{avg}}, f \geq f_{avg} \\ P_{c_{max}}, f < f_{avg} \end{cases} \tag{52}$$

Where:

$f_{max}$: the maximum fitness value in the population.

$f_{avg}$: the average fitness value in the population.

$\left[P_{c_{min}}, P_{c_{max}}\right]$: the specified range for crossover probability.

The chromosome crossover process employs the Partially Mapped Crossover (PMX) method. First, a random pair of chromosomes is selected, along with a set of start and end positions for several genes. It is essential to ensure that the selected positions on both chromosomes match. Then, two sets of genes are exchanged between the two chromosomes. Finally, conflict detection is performed on the chromosomes, and any conflicting genes are replaced until conflicts are resolved. The specific chromosome crossover process is illustrated in Fig. 7. This method helps to maintain the diversity of the population effectively.

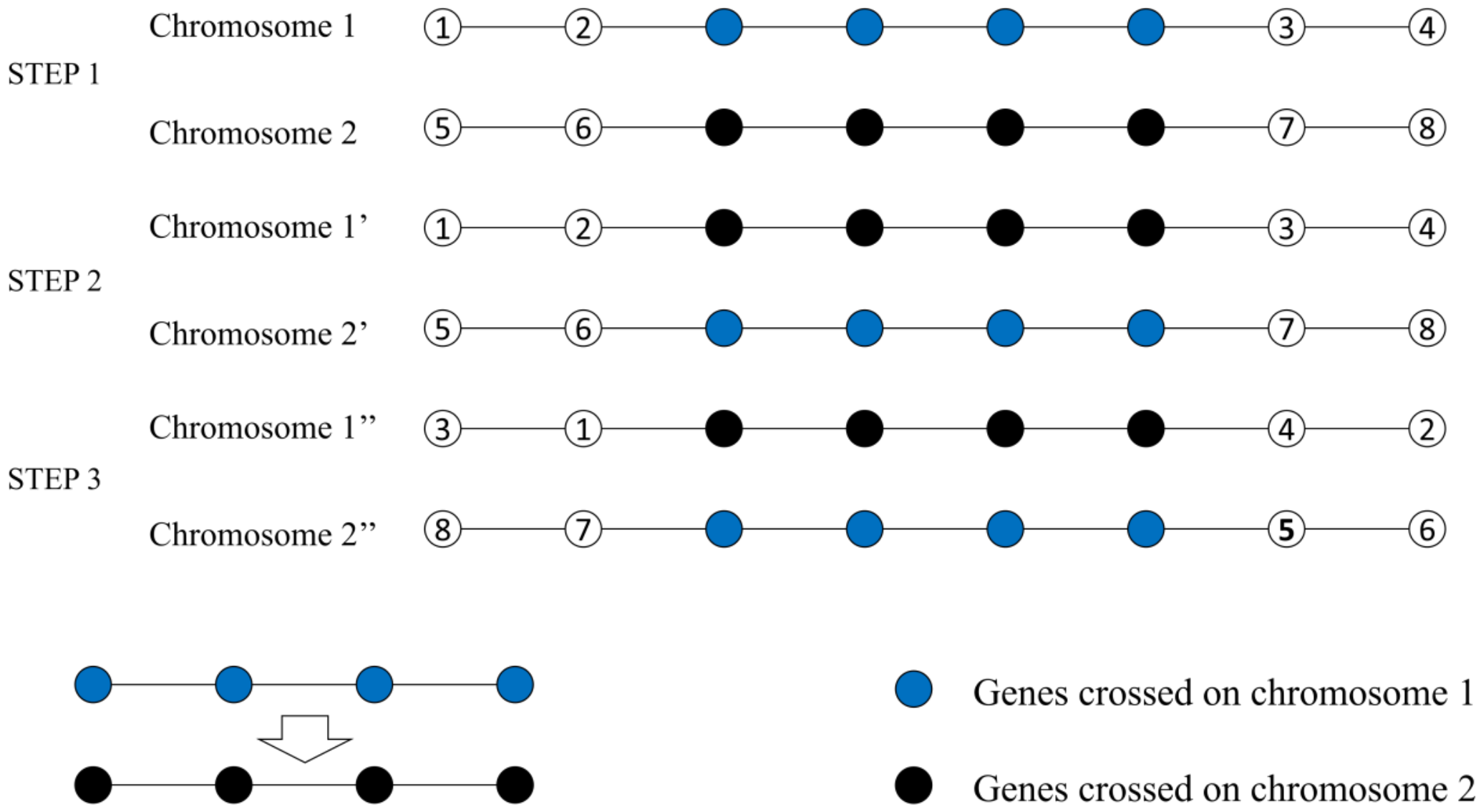

**Fig. 7.** Chromosome crossover process.

**c) Mutation**

In this paper, an adaptive mutation strategy is used to adjust the mutation probability, denoted as $P_m$, based on the fitness values, $f'_a$ and $f'_b$, of two individuals, $a$ and $b$, participating in the mutation operation. Here, $f'$ represents the maximum fitness value among the two individuals involved in the mutation operation, i.e., $f' = \max(f'_a, f'_b)$. When the population tends towards a local optimum, the value of $P_m$ will increase, and when the population is more dispersed, $P_m$ decreases. The expression for adjusting the adaptive mutation probability is:

$$P_m = \begin{cases} P_{m_{max}} - \frac{(P_{m_{max}} - P_{m_{min}})(f_{max} - f')}{f_{max} - f_{avg}}, f' \geq f_{avg} \\ P_{m_{max}}, f' < f_{avg} \end{cases} \tag{53}$$

Where $[P_{m_{min}}, P_{m_{max}}]$ represents the specified range for mutation probability.

The chromosome mutation process uses the strategy of gene segment reversal, which is superior to the traditional swap mutation strategy. This strategy can enhance the sorting performance and prevent feasible solutions from getting stuck in local optima too quickly. On chromosome 1, two crossover points, 1 and 9, are randomly selected. The gene segment between these two points is reversed to obtain the mutated chromosome 1', as shown in Fig. 8.

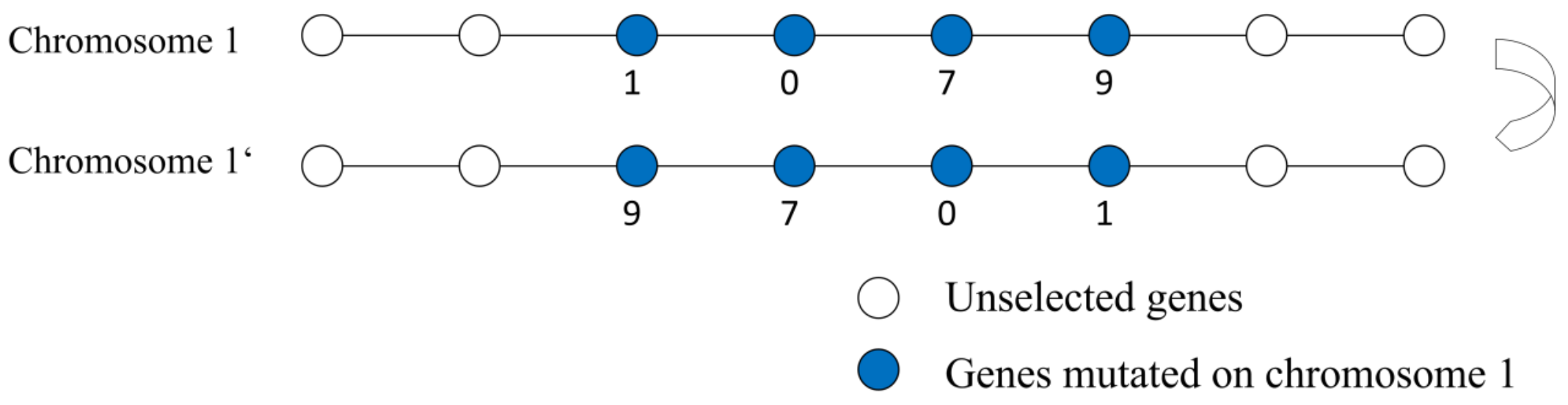

**Fig. 8.** Chromosome mutation process.

*5.2. Algorithm steps*

Building upon the detailed explanation of the genetic algorithm design in subsection 5.1, this subsection presents the specific operational steps of the improved genetic algorithm. The steps to perform the improved genetic algorithm are depicted in Fig. 9.

Step 1: Given the positions of parking lots, aircraft stands, and charging stations, calculate the distances between every pair of these locations. Encode the existing chromosomes.

Step 2: Generate an initial population based on a greedy approach. Set the iteration count as iter = 1.

Step 3: Calculate the fitness of each individual in the initial population using a fitness function and rank them. Apply an elite strategy by placing the best-performing individual into the elite set. Use a roulette wheel selection method to choose the remaining individuals for the next step.

Step 4: Determine the crossover probability, $P_c$, using an adaptive crossover strategy. Perform the crossover operation based on the partial-mapped crossover strategy. Individuals that meet the feasible solution conditions proceed to Step 6.

Step 5: Determine the mutation probability, $P_m$, using an adaptive mutation strategy. Perform the mutation operation based on the reversal mutation strategy. Individuals that meet the feasible solution conditions proceed to Step 6.

Step 6: Calculate the fitness of the new population and add the best-performing individual to the elite set.

Step 7: Check if the maximum iteration count has been reached. If not, increase the iteration count (iter = iter + 1) and return to Step 3. If the maximum iteration count has been reached, stop the algorithm and output the results.

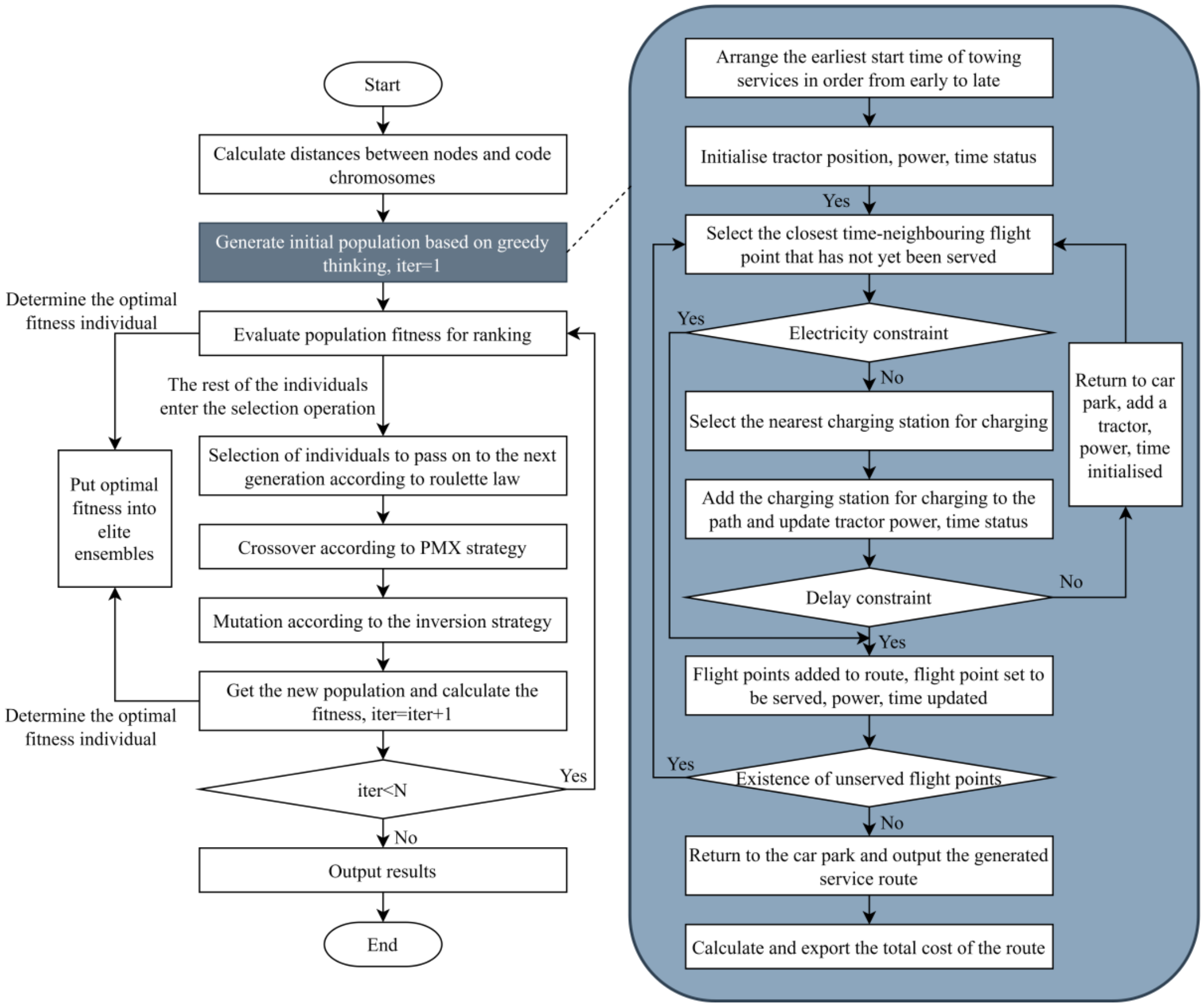

**Fig. 9.** Flowchart of the improved genetic algorithm.

## 6. Case Study

### *6.1. Case design*

Nanjing Lukou Airport stands as a prominent aviation hub in eastern China, receiving strong support from the Civil Aviation Administration of China. In 2022, it accommodated 12.14 million passengers, ranking 13th among mainland Chinese airports, and handled 378,000 tons of cargo, securing the 9th position in mainland China.

Nanjing Lukou Airport comprises two terminals, T1 and T2. T1, established in 1997, is primarily dedicated to domestic flights, whereas T2, inaugurated in 2014, serves international and regional flights. Given the independence of these two terminals and their significant differences in flight scale and configuration, two distinct scenarios have been designed to demonstrate the applicability of the methodology: scenario 1 (covering aircraft stands near T1 and some remote aircraft stands) and scenario 2 (covering aircraft stands near T2 and some remote aircraft stands).

**(1) Scenario Information**

Scenario 1 terminal features a finger-like terminal structure, with fewer close aircraft stands and more remote aircraft stands. The road network spans 4.21 kilometers and includes 18 close aircraft stands, 24 remote aircraft stand nodes, 3 charging station nodes, and 1 parking lot node. Vehicles have longer routing distances, and aircraft stands are relatively dispersed. Scenario 2 terminal features a front-row terminal structure, with more close aircraft stands and fewer remote aircraft stands. The road network extends 3.48 kilometers and comprises a total of 31 close aircraft stands, 8 remote aircraft stand nodes, 2 charging station nodes, and 1 parking lot node. Vehicles have shorter routing distances, and aircraft stands are more densely concentrated. The layout characteristics of the two scenarios are shown in Fig.10.

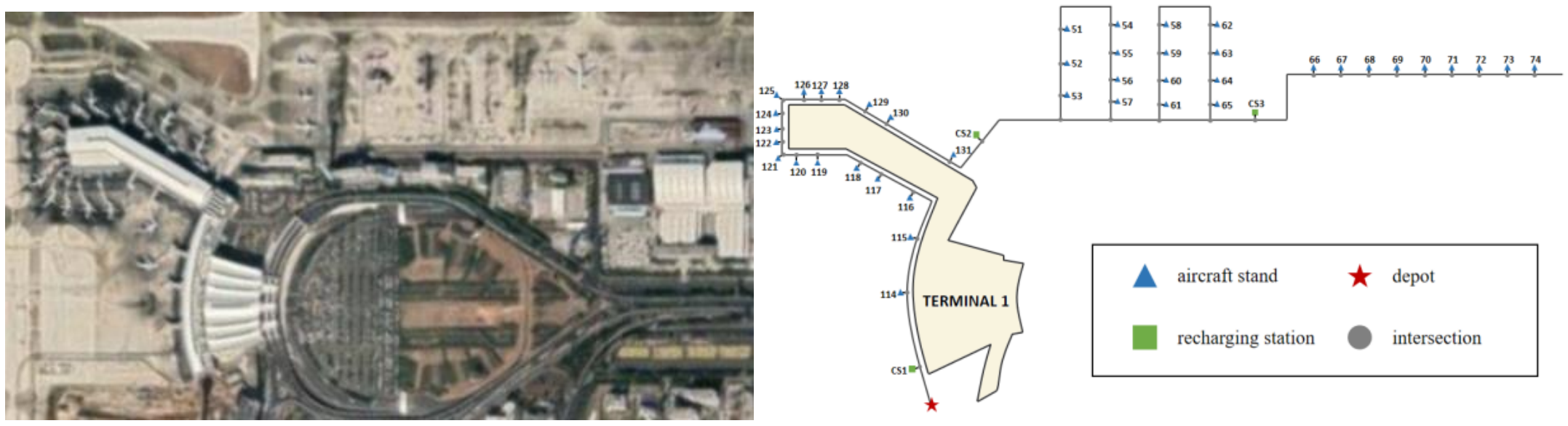

(a)Satellite view of scenario 1 (b) Network diagram of scenario 1

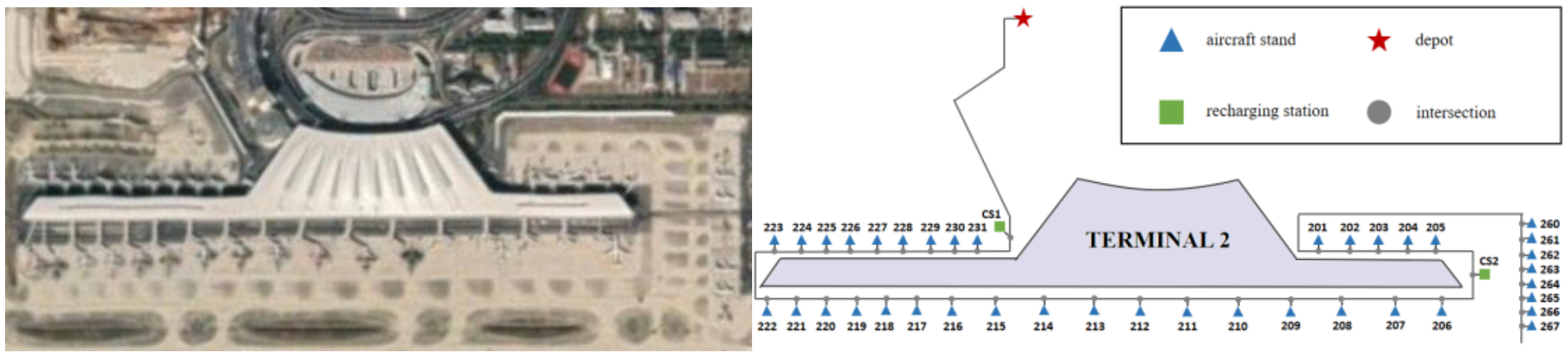

(c) Satellite view of scenario 2 (d) Network diagram of scenario 2

**Fig. 10.** The layout characteristics of the two scenarios.

### (2) Flight and tractor information

In scenario 1, there were a total of 103 transit flights on January 8, 2022, which requires a demand for 103 aircraft towing services. This averages approximately 8 services per hour. In scenario 2, there were 44 transit flights on January 8, 2022, requiring 44 aircraft towing services. This averages approximately 4 services per hour.

The flight input information for both scenarios includes the flight number, aircraft type (determining the aircraft mass), scheduled arrival time, scheduled departure time, and boarding gate. Detailed performance parameters for tractors and aircraft can be found in Appendix A.

## *6.2. Results*

### 6.2.1. *Algorithm comparison*

To validate the effectiveness of the proposed improved genetic algorithm, this section compares it with the traditional genetic algorithm. The program was run using Matlab R2018b. The parameters for the improved genetic algorithm are configured as follows: a population size of 200, a maximum iteration count of 1,000, a generation gap of 0.9, adaptive probability range for crossover in [0.6, 0.8], and adaptive probability range for mutation in [0.009, 0.2]. With a tractor speed of 25 km/h and a low-level charging strategy, both algorithms were run five times for scenarios 1 and 2 at Nanjing Lukou Airport. The resulting objective function values and solution time are compared in Table 2, and the iteration curves for the optimal solutions are depicted in Fig. 11 and Fig. 12.

**Table 2**

Comparison of objective function values and computation time obtained by both algorithms.

| Scenario number | Traditional genetic algorithm | | Improved genetic algorithm | |
|---|---|---|---|---|
| | Objective function value ($) | Computation time (s) | Objective function value ($) | Computation time (s) |
| T1-1 | 2190 | 291 | 1803 | 435 |
| T1-2 | 2608 | 258 | 1842 | 448 |
| T1-3 | 2233 | 279 | 1711 | 403 |
| T1-4 | 2326 | 283 | 1790 | 407 |
| T1-5 | 2270 | 277 | 1681 | 469 |
| Average | 2325 | 278 | 1765 | 432 |
| T2-1 | 643 | 136 | 546 | 244 |
| T2-2 | 665 | 138 | 535 | 240 |
| T2-3 | 576 | 128 | 550 | 241 |
| T2-4 | 622 | 132 | 548 | 239 |
| T2-5 | 592 | 127 | 538 | 234 |
| Average | 620 | 132 | 543 | 240 |

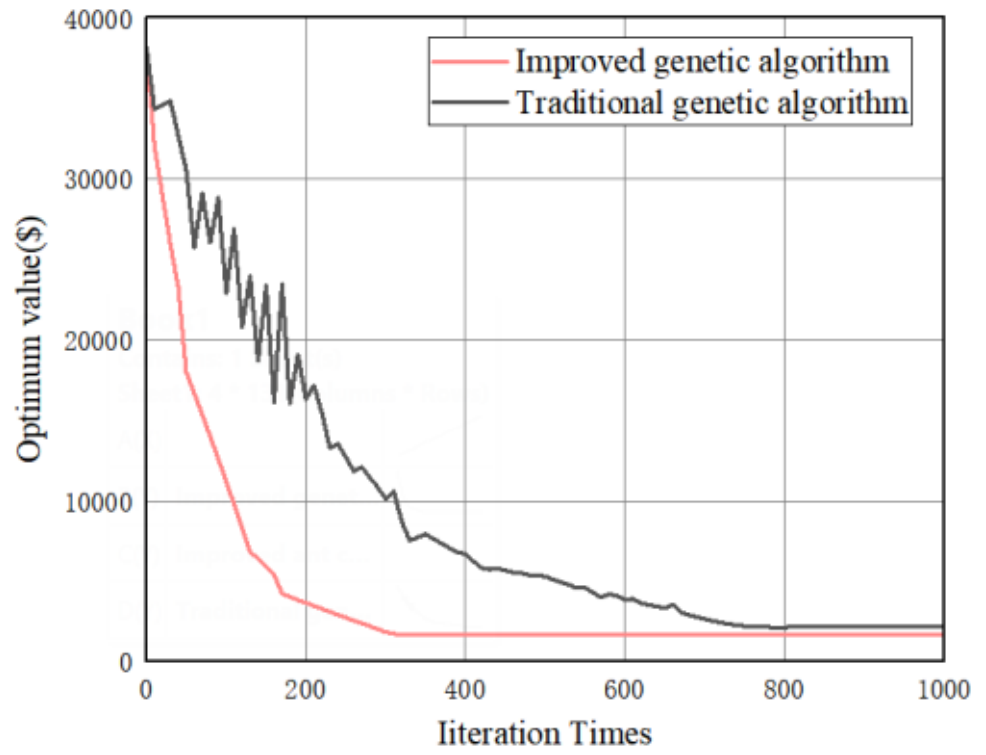

**Fig. 11.** Best performance curves for scenario 1.

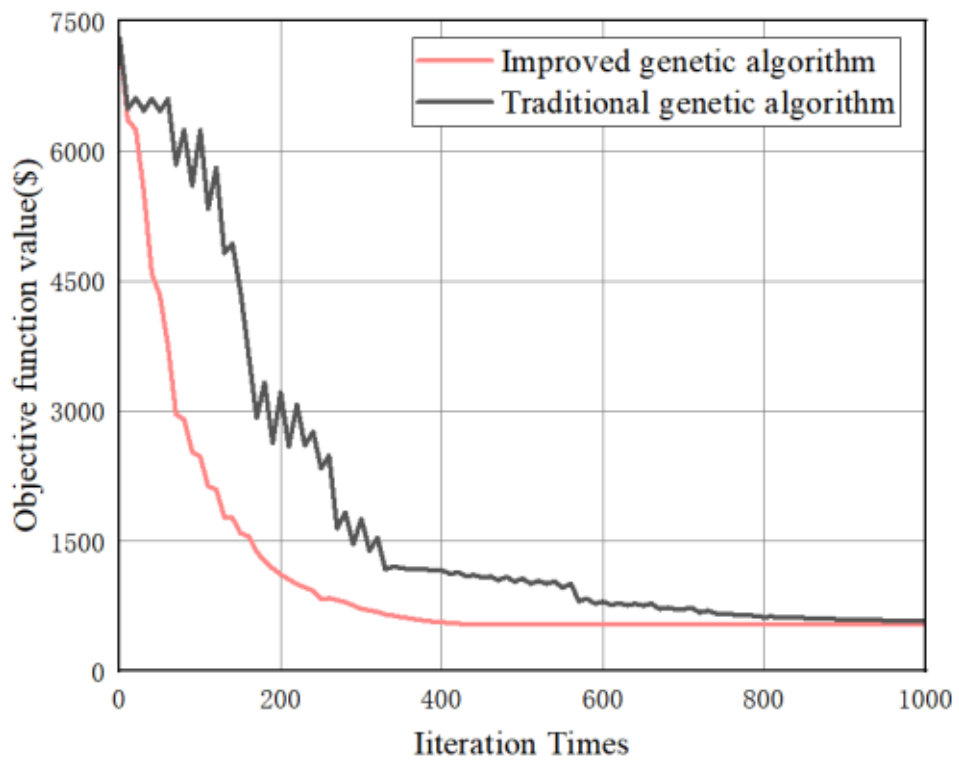

**Fig. 12.** Best performance curves for scenario 2.

Table 2 shows that, in scenario 1, the improved genetic algorithm results in a 24% reduction in the average objective function value when compared to the traditional genetic algorithm. For scenario 2, the use of the improved genetic algorithm leads to a 12% reduction in the average objective function value compared to the traditional genetic algorithm. In terms of computation times, in scenario 1, utilizing the improved genetic algorithm results in an average increase of 55% in computation time across the five runs. For scenario 2, the use of the improved genetic algorithm leads to an average increase of 82% in computation time across the five runs. By combining the information from Fig. 11 and Fig. 12, it can be concluded that the improved genetic algorithm converges more rapidly, and the results from each search are relatively consistent.

In summary, the improved genetic algorithm introduced in this paper, despite a minor drawback in computation time, notably enhances the algorithm's solution accuracy and convergence speed within an acceptable computation time extension. Therefore, we believe that this improved genetic algorithm is well-suited for solving the complex constraint-based airport electric tractor dispatching model presented in this paper.

### *6.2.2. Model comparison*

To validate the effectiveness of the proposed improved energy consumption model considering the start-stop process, this section compares it with the traditional energy consumption model. The same 25km/h pacing strategy and Low-level charging strategy are set, and the two models are used to solve scenario 1 and scenario 2 of Nanjing Lukou Airport five times, respectively. The objective function values and tractor fleet size comparisons are obtained as shown in Table 3.

**Table 3**

Comparison of objective function values and tractor fleet size obtained by both models.

| Scenario number | Traditional energy consumption model | | Improved energy consumption model | |
|---|---|---|---|---|
| | Objective function value ($) | Tractor Fleet size | Objective function value ($) | Tractor Fleet size |
| T1-1 | 1785 | 12 | 1803 | 12 |
| T1-2 | 1661 | 11 | 1842 | 13 |
| T1-3 | 1593 | 11 | 1711 | 12 |
| T1-4 | 1644 | 11 | 1790 | 13 |
| T1-5 | 1627 | 11 | 1681 | 12 |
| Average | 1662 | 11 | 1765 | 12 |
| T2-1 | 531 | 3 | 546 | 3 |
| T2-2 | 522 | 3 | 535 | 3 |
| T2-3 | 540 | 3 | 550 | 3 |
| T2-4 | 529 | 3 | 548 | 4 |
| T2-5 | 533 | 3 | 538 | 3 |
| Average | 531 | 3 | 543 | 3 |

Table 3 shows that, in scenario 1, the improved energy consumption model results in a 6% increase in the average objective function value when compared to the traditional energy consumption model. For scenario 2, the use of the improved energy consumption model leads to a 2% increase in the average objective function value compared to the traditional energy consumption model. In terms of tractor fleet size, utilizing the improved energy consumption model, the average required tractor fleet size for scenario 1 is 12 vehicles, which is an increase of 1 vehicle over the traditional energy consumption model. For scenario 2, there is little difference in the average required tractor fleet size between the two models, both being 3 vehicles.

In summary, the traditional energy consumption model can lead to a calculated fleet size that is insufficient to meet actual demands. There is a clear shortage of tractors, especially for scenario 1, which can result in the failure of tractor dispatching plans. During peak hours, the supply of tractors cannot keep up with the demand, ultimately leading to some flights not receiving timely service and causing delays.

### *6.2.3. Dispatching results*

Using the Electric Vehicle Routing Problem with Time Windows (EVRPTW) model established in this paper along with the improved genetic algorithm, assuming a 25 km/h pacing strategy and a low-level charging strategy. The optimal route plans for scenario 1 and scenario 2 at Nanjing Lukou Airport are obtained, as shown in Table 4.

**Table 4**

The optimal route plans for two scenarios.

| Instance | Tractor number | Route | Electricity consumption (kwh) | Travel distance (km) | Total cost ($) |
|---|---|---|---|---|---|
| Scenario 1 | E1 | d-15-23-24-27-39-42-51-62-63-70-cs2-68-80-48-85-100-101-102-cs1-103-d | 216.83 | 10.03 | 187.23 |
| | E2 | d-18-59-76-79-81-86-87-d | 82.76 | 2.90 | 137.59 |
| | E3 | d-1-7-14-33-35-13-43-44-46-cs2-71-72-d | 144.36 | 7.05 | 149.64 |
| | E4 | d-3-32-53-54-58-60-67-64-74-cs2-77-93-90-94-99-d | 169.14 | 8.26 | 181.25 |
| | E5 | d-17-11-28-40-38-d | 65.50 | 5.93 | 113.56 |
| | E6 | d-9-29-26-d | 37.51 | 2.57 | 88.05 |
| | E7 | d-16-55-61-65-d | 49.52 | 3.00 | 118.35 |
| | E8 | d-36-57-56-92-88-96-97-98-d | 94.00 | 3.05 | 146.03 |
| | E9 | d-19-22-41-69-91-89-d | 72.03 | 3.31 | 140.40 |
| | E10 | d-2-12-20-25-21-30-34-31-37-cs2-66-78-75-82-84-d | 171.60 | 9.88 | 173.90 |
| | E11 | d-6-49-50-45-52-47-73-83-95-d | 108.03 | 4.62 | 151.17 |
| | E12 | d-5-4-8-10-d | 48.96 | 2.85 | 76.56 |
| | Total | | 1260.24 | 63.45 | 1680.81 |
| Scenario 2 | E1 | d-1-7-8-2-22-24-23-26-27-cs2-29-32-33-34-35-36-d | 186.76 | 12.18 | 187.23 |
| | E2 | d-6-10-20-21-19-28-30-31-39-cs2-37-38-41-42-d | 164.64 | 5.83 | 165.44 |
| | E3 | d-4-3-5-11-9-12-13-14-15-16-cs1-18-17-25-40-43-44-d | 190.56 | 7.68 | 182.23 |
| | Total | | 541.96 | 25.69 | 534.90 |

According to the calculation results, scenario 1 requires a total of 12 electric tractors to provide towing services. Tractor 1 visits the charging station twice during its journey, while tractors 3, 4, and 10 visit the charging station once during their journeys. The remaining tractors do not visit the charging station during their journeys. The total distance covered to complete all towing services is 63.45 kilometers, the total energy consumption is 1260.24 kilowatt-hours, and the total cost is 1680.81 USD. For scenario 2, a total of three electric tractors are needed to provide towing services, and all three tractors visit the charging station once during their journeys. The total distance covered to complete all towing services is 25.69 kilometers, the total energy consumption is 541.96 kilowatt-hours, and the total cost is 534.90 USD. The travelling paths of the tractors in these two scenarios are visualized in Fig. 13(a) and (b).(Take E10 and E2 for example).

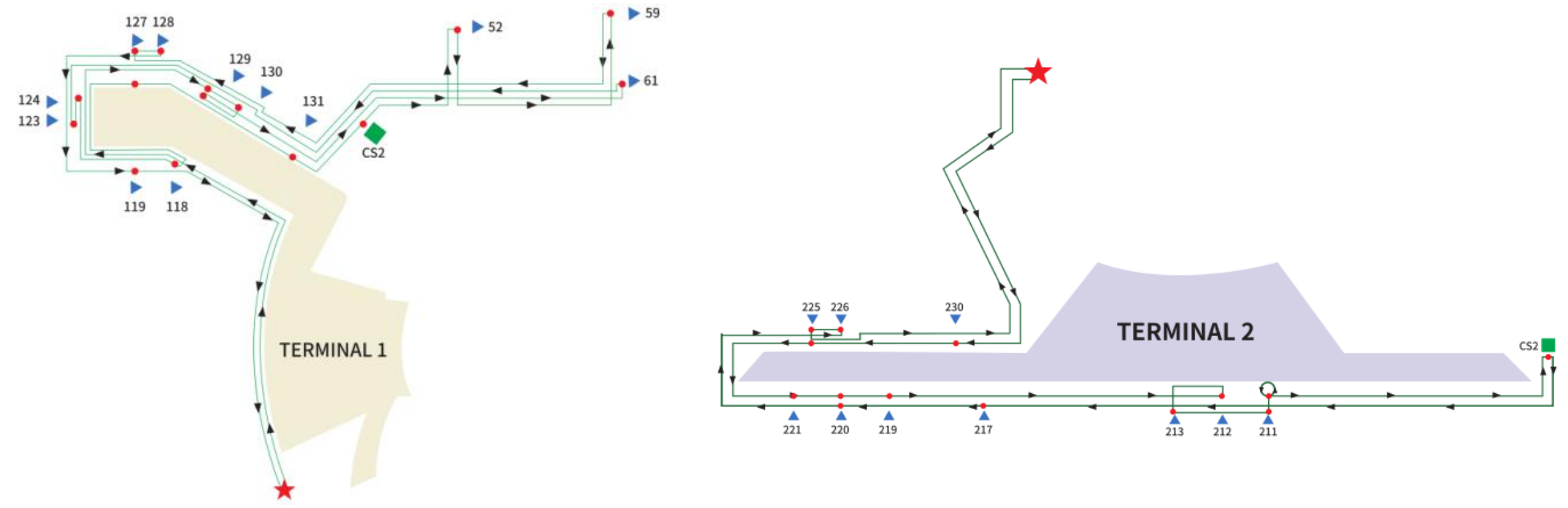

(a) (b)

**Fig. 13.** Viewable travelling paths a) for E10 in scenario 1; b) for E2 in scenario 2.

# 7. Discussion

## *7.1. Charge-discharge coupling strategy for dispatching in scenario 1*

This section discusses the most suitable pacing strategy for three charging levels: high (γ =1), medium (γ =0.9), and low (γ =0.8) in scenario 1, ultimately determining the optimal coupling strategy for charging and discharging. The results for different charge-discharge coupling strategies in scenario 1 are summarized in Table 5.

**Table 5**

Computation results of different charge-discharge coupling strategies for scenario 1.

| Charging strategy | Speed (km/h) | Required tractors | Fixed cost ($) | Variable cost ($) | Charging Cost ($) | Time cost ($) | Total cost ($) |
|---|---|---|---|---|---|---|---|
| High-level (γ =1) | 5 | 13 | 650.00 | 329.27 | 252.08 | 510.85 | 1742.20 |
| | 10 | 13 | 650.00 | 332.98 | 252.55 | 540.57 | 1776.10 |
| | 15 | 13 | 650.00 | 335.03 | 252.77 | 595.78 | 1833.58 |
| | 20 | 13 | 650.00 | 348.98 | 253.75 | 620.57 | 1873.30 |
| | 25 | 12 | 600.00 | 337.73 | 253.41 | 489.24 | 1680.38 |
| Mid-level (γ =0.9) | 5 | 13 | 650.00 | 325.51 | 251.90 | 501.65 | 1729.06 |
| | 10 | 13 | 650.00 | 331.59 | 252.30 | 510.16 | 1744.05 |
| | 15 | 13 | 650.00 | 326.74 | 252.17 | 562.72 | 1791.63 |
| | 20 | 13 | 650.00 | 341.69 | 253.37 | 588.97 | 1834.03 |
| | 25 | 12 | 600.00 | 309.24 | 251.53 | 492.51 | 1653.28 |
| Low-level (γ =0.8) | 5 | 13 | 650.00 | 345.99 | 253.18 | 574.81 | 1823.98 |
| | 10 | 13 | 650.00 | 332.23 | 252.44 | 602.16 | 1836.83 |
| | 15 | 13 | 650.00 | 327.38 | 252.22 | 639.61 | 1869.21 |
| | 20 | 12 | 600.00 | 318.23 | 252.21 | 492.16 | 1663.60 |
| | 25 | 12 | 600.00 | 317.25 | 252.05 | 511.51 | 1680.81 |

From Table 5, it shows that when the charging strategy is set to the high level with γ = 1, as the speed increases from 5 km/h to 20 km/h, the total number of required tractors remains at 13, and the total cost significantly rises from $1742.20 to $1873.30. At a speed of 25 km/h, the 13 tractors surpass the service demand, and the fleet can be adjusted to 12. The reduction in the number of required tractors leads to a sharp decrease in the total cost, reaching a minimum of $1680.38.

When the charging strategy is at the medium level with γ = 0.9, a turning point occurs as the speed increases from 20 km/h to 25 km/h. Before this turning point, the total number of required tractors remains at 13, and the total cost shows a slight upward trend as the speed increases, going from \$1729.06 to \$1834.03. After the turning point, the number of vehicles is adjusted from 13 to 12, and the total cost reaches a minimum of \$1653.28.

For the low charging strategy with γ = 0.8, the turning point occurs as the speed increases from 15 km/h to 20 km/h. Similarly, the total number of required tractors is adjusted from 13 to 12, but the turning point occurs at a different speed. Moreover, there are differences in the pacing strategy that results in the minimum total cost. In this case, a speed of 20 km/h is optimal, resulting in a total cost of \$1663.60. However, with a subsequent increase to 25 km/h, the total cost slightly rises to \$1680.81.

In summary, the total cost variation under different charge-discharge coupling strategies in scenario 1 is depicted in Fig. 14.

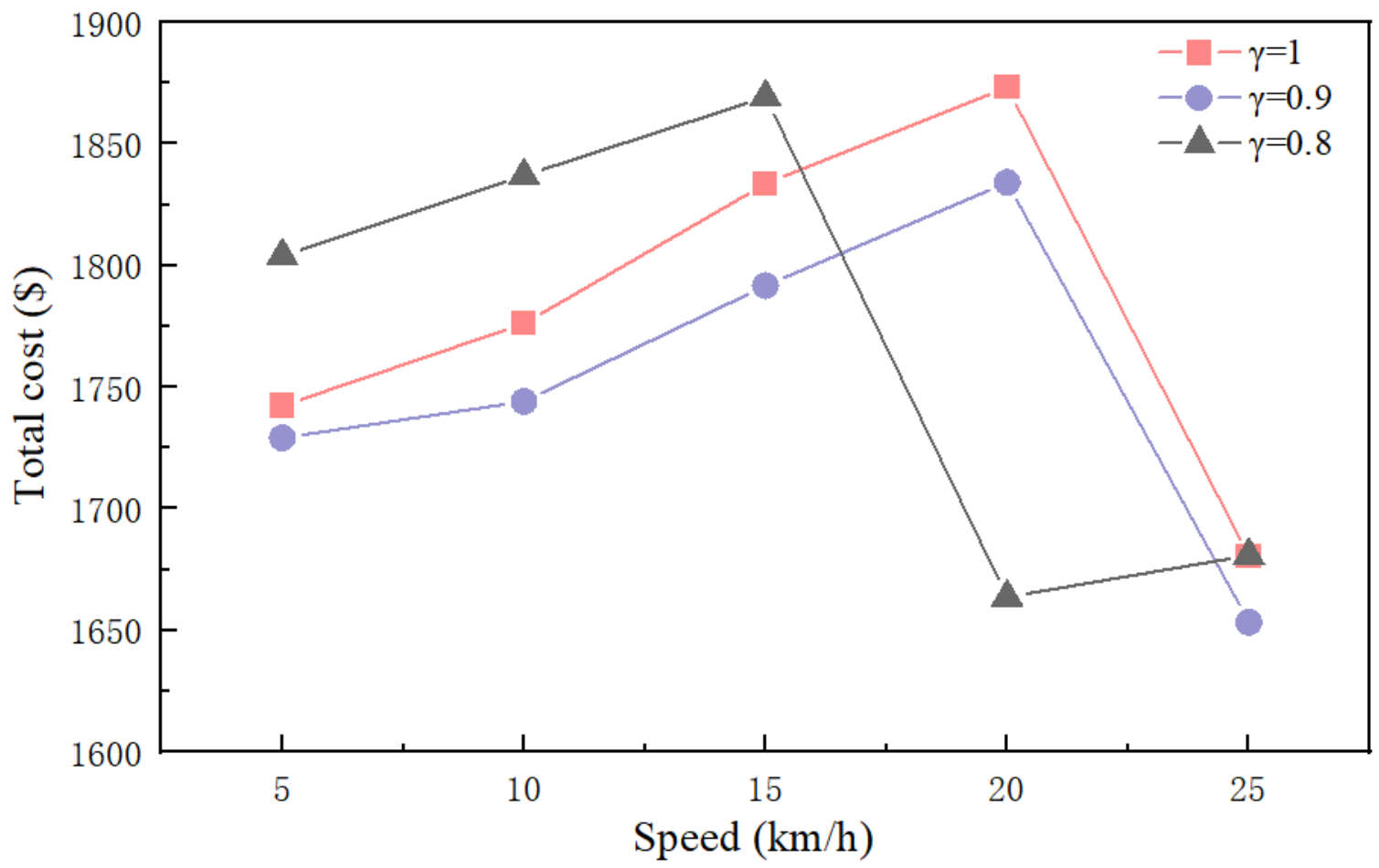

**Fig. 14.** The total cost under differential charge-discharge coupling strategies in scenario 1.

Fig. 14 clearly indicates that, when maintaining a consistent speed, the mid-level charging strategy generally outperforms other strategies. When evaluating all charge-discharge coupling strategies, it shows that a speed of 25 km/h combined with the mid-level charging strategy is the optimal choice. This configuration yields the lowest total cost, reducing it to \$1653.28, and involves a total of 12 tractors.

*7.2. Charge-discharge coupling strategy for dispatching in scenario 2*

This section discusses the suitable pacing strategy for three charging levels: high (γ =1), medium (γ =0.9), and low (γ =0.8) in scenario 2, and determining the optimal coupling strategy for charging and discharging. The results for different charge-discharge coupling strategies in scenario 2 are summarized in Table 6.

**Table 6**

Computation results of different charge-discharge coupling strategies for scenario 2.

| Charging strategy | Speed (km/h) | Required tractors | Fixed cost ($) | Variable cost ($) | Charging Cost ($) | Time cost ($) | Total cost ($) |
|---|---|---|---|---|---|---|---|
| High-level (γ =1) | 5 | 5 | 250.00 | 144.39 | 109.10 | 158.42 | 661.91 |
| | 10 | 5 | 250.00 | 141.64 | 108.97 | 177.59 | 678.20 |
| | 15 | 5 | 250.00 | 152.85 | 109.77 | 168.68 | 681.30 |
| | 20 | 5 | 250.00 | 140.63 | 109.06 | 189.71 | 689.40 |
| | 25 | 4 | 200.00 | 136.15 | 108.90 | 120.70 | 565.75 |
| Mid-level (γ =0.9) | 5 | 4 | 200.00 | 136.61 | 108.55 | 140.61 | 584.95 |
| | 10 | 4 | 200.00 | 134.84 | 108.51 | 146.39 | 589.74 |
| | 15 | 4 | 200.00 | 147.94 | 109.46 | 176.98 | 634.38 |
| | 20 | 3 | 150.00 | 136.27 | 108.78 | 141.53 | 536.58 |
| | 25 | 3 | 150.00 | 138.53 | 109.22 | 152.46 | 550.21 |
| Low-level (γ =0.8) | 5 | 4 | 200.00 | 142.96 | 109.00 | 137.61 | 589.57 |
| | 10 | 4 | 200.00 | 127.48 | 108.02 | 167.82 | 603.32 |
| | 15 | 3 | 150.00 | 131.96 | 108.33 | 120.20 | 510.48 |
| | 20 | 3 | 150.00 | 128.98 | 108.51 | 140.90 | 528.39 |
| | 25 | 3 | 150.00 | 128.45 | 108.39 | 148.06 | 534.90 |

From Table 6, it shows that when the charging strategy is set to the high level with $\gamma = 1$, as the speed increases from 5 km/h to 20 km/h, the total number of required tractors remains at 5. The total cost slightly increases, going from \$661.91 to \$689.40. A turning point occurs when the speed reaches 25 km/h, with 5 tractors exceeding the service demand and can be adjusted to 4. At this point, the total cost reaches its minimum at \$565.75.

When the charging strategy is at the medium level with $\gamma = 0.9$, a turning point similarly occurs, but the speed corresponding to this turning point is different, happening when the speed increases from 15 km/h to 20 km/h. Using the same number of tractors, all indicators show a consistent upward trend with increasing speed. Moreover, the optimal speed resulting in the minimum total cost differs. When $\gamma = 0.9$, a speed of 20 km/h is more suitable, resulting in a total cost of \$536.58.

With the charging strategy set to the low level and $\gamma = 0.8$, a turning point occurs as the speed increases from 10 km/h to 15 km/h. Meantime, the total number of required tractors decreases from 4 to 3, and the total cost decreases from \$603.32 to \$510.48. Subsequently, at 15 km/h, 20 km/h, and 25 km/h, the total number of required tractors remains at 3, and the total cost shows a slight upward trend, rising to \$534.90. In summary, the total cost under different charge-discharge coupling strategies in scenario 2 is illustrated in Fig. 15.

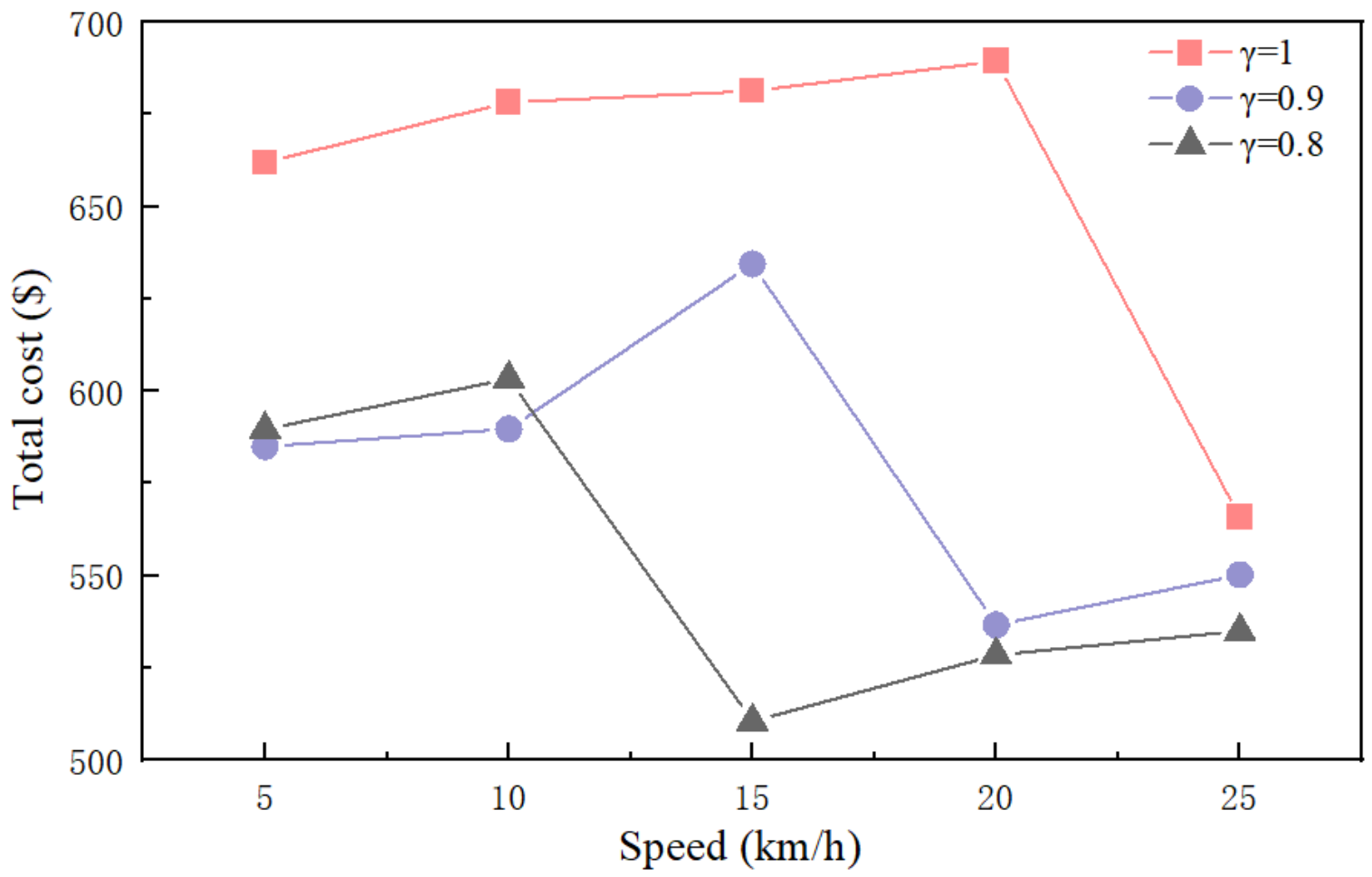

**Fig. 15.** The total cost under differential charge-discharge coupling strategies in scenario 2.

From Fig. 15, it shows that when the speed remains consistent, the low-level charging strategy is notably superior to other strategies. When comparing all charge-discharge coupling strategies, it is found that a speed of 15 km/h combined with the low-level charging strategy is the most optimal choice. This configuration leads to the lowest total cost, reducing it to $510.48, and necessitating a total of three tractors.

In summary, for scenario 1, a speed of 25 km/h with the mid-level charging strategy may be considered, while for scenario 2, a speed of 15 km/h with the low-level charging strategy is advisable. Appropriately slowing down the driving speed of electric tractors during towing service tasks is advantageous for cost reduction and optimizing tractor utilization. Nevertheless, an overly drastic reduction in speed might result in tractors being unable to complete their tasks within the designated time, thereby requiring an increase in the fleet size.

### *7.3. Impact of future flight scale changes on strategies*

In Sections 7.1 and 7.2, we have explored the appropriate charge-discharge coupling strategy for the current scenarios of both terminals at Nanjing Lukou Airport using real data. However, global economic growth and advancements in aviation technology have resulted in an annual increase in flight scales, particularly due to the gradual lifting of pandemic restrictions worldwide, leading to a significant rise in international flight scales at Terminal 2. Additionally, as Terminal 1's capacity gradually reaches saturation, some domestic flights will be transferred to Terminal 2. These factors will lead to significant changes in the flight scales at airport terminals in the future. Therefore, further research is needed to examine the applicability of current strategies considering future changes in flight scales and their potential impacts.

In this section, we predict changes in flight scales over the next three and five years based on the annual growth rate of aircraft movements at the airport and terminal capacity constraints (detailed forecasting process can be found in Appendix B). Specifically, we categorize the scenarios based on two distinct timeframes: the near-term future case, encompassing the upcoming three years, and the mid-term future case, extending to the upcoming five years. In scenario 1, the flight scale for the near-term future case is 130, and for the mid-term future case, it is 150. In scenario 2, the flight scale for the near-term future case is 80, and for the mid-term future case, it is 130. Then, we randomly determine turnaround times based on historical patterns of peak and off-peak flight periods at Nanjing Lukou Airport. In addition, considering the existing aircraft stands at the airport, we randomly allocate an aircraft stand for each flight (note that each aircraft stand can only accommodate one aircraft at any given time). Finally, we derive the optimal charge-discharge coupling strategies for the two scenarios as shown in Tables 7 and 8.

**Table 7**

Optimal charge-discharge coupling strategies for scenario 1 in the future.

| Charging strategy | Speed (km/h) | Current (103) | | Near-term future case (130) | | Mid-term future case (150) | |
|---|---|---|---|---|---|---|---|
| | | Required tractors | Total cost ($) | Required tractors | Total cost ($) | Required tractors | Total cost ($) |
| High-level (γ =1) | 5 | 13 | 1742.20 | 16 | 2311.62 | 20 | 2913.92 |
| | 10 | 13 | 1776.10 | 16 | 2363.93 | 20 | 2942.52 |
| | 15 | 13 | 1833.58 | 16 | 2383.64 | 20 | 2936.12 |
| | 20 | 13 | 1873.30 | 16 | 2420.76 | 20 | 2974.21 |
| | 25 | 12 | 1680.38 | 15 | 2228.62 | 19 | 2805.31 |
| Mid-level (γ =0.9) | 5 | 13 | 1729.06 | 16 | 2298.63 | 20 | 2927.41 |
| | 10 | 13 | 1744.05 | 16 | 2339.37 | 20 | 2943.37 |
| | 15 | 13 | 1791.63 | 16 | 2353.13 | 20 | 2976.39 |
| | 20 | 13 | 1834.03 | 16 | 2386.35 | 20 | 3012.64 |
| | 25 | 12 | 1653.28 | 15 | 2207.14 | 19 | 2829.66 |
| Low-level (γ =0.8) | 5 | 13 | 1823.98 | 16 | 2391.64 | 20 | 2979.73 |
| | 10 | 13 | 1836.83 | 16 | 2419.87 | 20 | 3004.85 |
| | 15 | 13 | 1869.21 | 16 | 2452.68 | 20 | 3022.67 |
| | 20 | 12 | 1663.60 | 15 | 2232.52 | 20 | 3040.73 |
| | 25 | 12 | 1680.81 | 15 | 2281.93 | 19 | 2856.94 |

**Table 8**

Optimal charge-discharge coupling strategies for scenario 2 in the future.

| Charging strategy | Speed (km/h) | Current (44) | | Near-term future case (80) | | Mid-term future case (130) | |
|---|---|---|---|---|---|---|---|
| | | Required tractors | Total cost ($) | Required tractors | Total cost ($) | Required tractors | Total cost ($) |
| | 5 | 5 | 661.91 | 9 | 1191.67 | 16 | 2195.21 |
| High-level (γ =1) | 10 | 5 | 678.20 | 9 | 1218.92 | 16 | 2247.84 |
| | 15 | 5 | 681.30 | 8 | 1091.22 | 16 | 2305.63 |
| | 20 | 5 | 689.40 | 8 | 1136.71 | 15 | 2065.61 |
| | 25 | 4 | 565.75 | 8 | 1143.37 | 15 | 2099.98 |
| Mid-level (γ =0.9) | 5 | 4 | 584.95 | 9 | 1205.09 | 16 | 2232.32 |
| | 10 | 4 | 589.74 | 8 | 1066.95 | 16 | 2285.72 |
| | 15 | 4 | 634.38 | 8 | 1092.79 | 15 | 2041.77 |
| | 20 | 3 | 536.58 | 8 | 1094.55 | 15 | 2079.14 |
| | 25 | 3 | 550.21 | 8 | 1126.87 | 15 | 2113.56 |
| Low-level (γ =0.8) | 5 | 4 | 589.57 | 9 | 1217.88 | 16 | 2255.73 |
| | 10 | 4 | 603.32 | 8 | 1073.30 | 16 | 2331.58 |
| | 15 | 3 | 510.48 | 8 | 1090.33 | 15 | 2088.45 |
| | 20 | 3 | 528.39 | 8 | 1114.22 | 15 | 2106.47 |
| | 25 | 3 | 534.90 | 8 | 1135.11 | 15 | 2147.89 |

According to Table 7, for scenario 1, the near-term future case is best suited for a 25 km/h speed and mid-level charging strategy. For this context, the total cost reaches a minimum of $2,207.14, with a corresponding total of 15 tractors. For the mid-term future case, a 25 km/h speed and high-level charging strategy are recommended, resulting in a total cost of $2,805.31 and 19 tractors.

According to Table 8, for scenario 2, the near-term future case is optimal with a 10 km/h speed and mid-level charging strategy. Here, the total cost reaches a minimum of $1,066.95, with a total of 8 tractors. As for the mid-term future case, a 15 km/h speed and mid-level charging strategy are advised, resulting in a total cost of $2,041.77 and 15 tractors.

In conclusion, the strategies are not set in stone. As future flight scales increase, both scenarios can moderately increase their charging strategies based on the current ones to achieve optimal benefits. Additionally, it is worth noting that although the flight scale for the near-term future case in scenario 1 is the same as the flight scale for the mid-term future case in scenario 2, there are differences in the optimal charge-discharge coupling strategies between these two scenarios. Specifically, scenario 1 is best suited for a 25 km/h speed and mid-level charging strategy, while scenario 2 is best suited for a 15 km/h speed and mid-level charging strategy. This indicates that the configuration of aircraft stands has a certain influence on the design of charge-discharge coupling strategies. Compared to airports with densely concentrated aircraft stands, airports with dispersed aircraft stand are better suited for adopting higher charging rates.

## 8. Conclusion

This paper has introduced a method for dispatching electric ground service tractors at airports, taking into account complex energy consumption and piecewise linear charging and has evaluated the impact of different charge-discharge coupling strategies on operational costs.

First, we have established a nonlinear energy consumption model that considers the start-stop process of electric tractors and their service status. We have also created an airport electric tractor dispatching model with time windows based on the piecewise linear charging function. The models have been proven to more accurately reflect the discharging and charging processes of electric tractors, aligning more closely with real-world scenarios.

Second, we have designed an improved genetic algorithm to solve the models, incorporating greedy thinking and adaptive strategies, ultimately obtaining the optimal tractor dispatching solution. To validate the effectiveness of the algorithm, we have compared it with the traditional genetic algorithm. The results have shown that the improved genetic algorithm significantly enhances solution accuracy and convergence speed.

Finally, we have applied our method to two scenarios at Nanjing Lukou Airport (T1 and T2), developing differentiated charge-discharge coupling strategies. Additionally, we have explored potential fluctuations in flight scales in the future, providing insights for the long-term development of optimal strategies. Our research has revealed the following:

1. As the size of future flights increases, airports can appropriately increase the charging strategy.
2. Compared to airports with dense aircraft stands, airports with dispersed aircraft stands are suitable for higher pacing strategies.
3. Reducing driving speed offers cost-saving advantages, but an overly aggressive pursuit of lower speeds may result in an increase in fleet size and total costs.

Our proposed method significantly enhances model accuracy and offers theoretical guidance for the development of differentiated charge-discharge coupling strategies for electric tractors at airports. However, the study has not extended the application of this method to other types of ground service vehicles. Further research can easily adjust the models and procedures to confirm the relevance of these conclusions for different vehicles involved in airport ground service processes.

## Appendix A

Table A.1 shows performance parameters of tractors.

**Table A.1**

Performance parameters of tractors.

| Parameters | Meaning | Numerical value |
| --- | --- | --- |
| $w$ | Distance of towing service | 50m |
| $u$ | Maximum allowable towing speed | 10km/h |
| $a$ | Magnitude of acceleration during the travel process | 0.8m/s$^2$ |
| $\varepsilon$ | Magnitude of acceleration during the towing process | 0.6m/s$^2$ |
| $g$ | Gravitational acceleration | 9.81N/kg |
| $\rho$ | Air density | 1.2041kg/m$^3$ |
| $A$ | Frontal area of the towing tractor | 3.912m$^2$ |
| $B$ | Battery capacity | 150kw·h |
| $C_d$ | Aerodynamic towing coefficient | 0.7 |
| $C_r$ | Rolling friction coefficient | 0.03 |
| $m$ | Tractor mass | 15300kg |
| $\emptyset^d$ | Motor discharge output efficiency | 1.184692 |
| $\emptyset^r$ | Motor recovery power input efficiency | 0.846055 |
| $\varphi^d$ | Battery energy output efficiency | 1.112434 |
| $\varphi^r$ | Battery energy recovery efficiency | 0.928465 |
| $c_1$ | Fixed cost per tractor | 50$/tractor |
| $c_2$ | Charging cost per unit kilowatt-hour | 0.2$/kW·h |
| $c_3$ | Maintenance cost per unit meter | 0.005$/m |
| $c_e$ | Waiting cost per unit minute | 0.1$/min |
| $c_l$ | Delay cost per unit minute | 0.5$/min |

Table A.2 shows performance parameters of aircraft.

**Table A.2**

Performance parameters of aircraft.

| Parameters | Medium aircraft | Heavy aircraft | Super heavy aircraft |
|---|---|---|---|
| $q_{AIR}$ | 175kw | 300kw | 350kw |
| $t_{AIR}$ | 2min | 2min | 2min |
| $q_{ILL}$ | 3.7kw | 14.1kw | 33.9kw |
| $t_{ILL}$ | 2min | 2min | 2min |
| $q_{LAU}$ | 384kw | 783kw | 783kw |
| $t_{LAU}$ | 0.75min | 0.75min | 0.75min |

(Medium aircraft: weight≤A321; Heavy aircraft: A321 < weight < B747; Super heavy aircraft: weight≥B747)

## Appendix B

This section provides a detailed overview of the process of forecasting short-term and mid-term future flight scales for two scenarios at Nanjing Lukou Airport.

**(1) Establishing baseline data**

Firstly, we gather the annual growth rate of takeoffs and landings at Nanjing Lukou Airport before the pandemic. According to the Civil Aviation Administration's annual statistics report, the average annual growth rate of takeoffs and landings from 2014 to 2018 is 14%.

Secondly, we compute the flight recovery rate post-pandemic for Nanjing Lukou Airport. We collect the latest data changes in flight scale for two scenarios: the annual recovery rate for scenario 1 is approximately 12%, and for scenario 2, it is approximately 30%. Given that the flight recovery rate decreases annually, we assume a decrease rate of 6%.

Finally, we assume that one year later, scenario 1 will reach saturation. With improvements in service quality and technology, the saturated terminal can appropriately increase the flight scale, resulting in a growth rate of approximately 7% after saturation. The remaining flights will be transferred to Terminal 2.

**(2) Establishing growth models**

We employ growth model to forecast future flight scales. At the onset of the pandemic relaxation, the flight scales for scenario 1 and scenario 2 are $T_{1cur} = 103$ and $T_{2cur} = 44$ , respectively. We denote the growth rates with the following symbols:

$r_{norm}$: Annual pre-pandemic flight growth rate.

$r_{1rec}$: Initial flight recovery rate for scenario 1 post-pandemic relaxation.

$r_{2rec}$: Initial flight recovery rate for scenario 2 post-pandemic relaxation.

$r$: Annual decrease rate in flight recovery rate post-pandemic relaxation.

$r_{sat}$: Annual flight growth rate for the saturated terminal.

**(3) Deriving flight scale for short-term future case**

For scenario 1:

$$\boldsymbol{T_{1short} = T_{1cur} * (1 + r_{1rec})(1 + r_{sat})(1 + r_{sat}) = 132.1} \quad \text{(B1)}$$

For scenario 2:

$$\boldsymbol{T_{2short} = T_{2cur} * (1 + r_{2rec})(1 + r_{2rec} - r)(1 + r_{2rec} - 2r) = 83.7} \quad \text{(B2)}$$

**(4) Deriving flight scale for mid-term future case**

For scenario 1:

$$\boldsymbol{T_{1mid} = T_{1short} * (1 + r_{sat})^2 = 151.2} \quad \text{(B3)}$$

For scenario 2:

$$\boldsymbol{T_{2mid} = T_{2short} * (1 + r_{norm})^2 + T_{1short} * (1 + r_{norm})^2 - T_{1mid} = 129.3} \quad \text{(B4)}$$

For the purpose of simulating data, we round the predicted results to the nearest ten. Consequently, in scenario 1, the short-term future flight scale is 130, and the mid-term future flight scale is 150. In scenario 2, the short-term future flight scale is 80, and the mid-term future flight scale is 130.

## References

Adler, Jonathan D., and Pitu B. Mirchandani. 2014. "Online Routing and Battery Reservations for Electric Vehicles with Swappable Batteries." *Transportation Research Part B: Methodological* 70:285–302. doi: 10.1016/j.trb.2014.09.005.

Amini, M. Hadi, Amin Kargarian, and Orkun Karabasoglu. 2016. "ARIMA-Based Decoupled Time Series Forecasting of Electric Vehicle Charging Demand for Stochastic Power System Operation." *Electric Power Systems Research* 140:378–90. doi: 10.1016/j.epsr.2016.06.003.

Bao, Dan-Wen, Jia-Yi Zhou, Zi-Qian Zhang, Zhuo Chen, and Di Kang. 2023. "Mixed Fleet Scheduling Method for Airport Ground Service Vehicles under the Trend of Electrification." *Journal of Air Transport Management* 108:102379. doi: 10.1016/j.jairtraman.2023.102379.

Basso, Rafael, Balázs Kulcsár, Bo Egardt, Peter Lindroth, and Ivan Sanchez-Diaz. 2019. "Energy Consumption Estimation Integrated into the Electric Vehicle Routing Problem." *Transportation Research Part D: Transport and Environment* 69:141–67. doi: 10.1016/j.trd.2019.01.006.

Bektaş, Tolga, and Gilbert Laporte. 2011. "The Pollution-Routing Problem." *Transportation Research Part B: Methodological* 45(8):1232–50. doi: 10.1016/j.trb.2011.02.004.

Brevoord, Julia M. 2021. "Electric Vehicle Routing Problems: The Operations of Electric Towing Trucks at an Airport under Uncertain Arrivals and Departures." Retrieved June 7, 2023 (https://essay.utwente.nl/86091/).

Conrad, Ryan, and Miguel Figliozzi. 2011. "The Recharging Vehicle Routing Problem." *Proc. of the 61st Annual IIE Conference*.

Du, Jia Yan, Jens O. Brunner, and Rainer Kolisch. 2014. "Planning Towing Processes at Airports More Efficiently." *Transportation Research Part E: Logistics and Transportation Review* 70:293–304. doi: 10.1016/j.tre.2014.07.008.

Felipe, Ángel, M. Teresa Ortuño, Giovanni Righini, and Gregorio Tirado. 2014. "A Heuristic Approach for the Green Vehicle Routing Problem with Multiple Technologies and Partial Recharges." *Transportation Research Part E: Logistics and Transportation Review* 71:111–28. doi: 10.1016/j.tre.2014.09.003.

Goeke, Dominik, and Michael Schneider. 2015. "Routing a Mixed Fleet of Electric and Conventional Vehicles." *European Journal of Operational Research* 245(1):81–99. doi: 10.1016/j.ejor.2015.01.049.

Guo, Weian, Ping Xu, Zhen Zhao, Lei Wang, Lei Zhu, and Qidi Wu. 2020. "Scheduling for Airport Baggage Transport Vehicles Based on Diversity Enhancement Genetic Algorithm." *Natural Computing* 19(4):663–72. doi: 10.1007/s11047-018-9703-0.

Han, Xue, Peixin Zhao, and Dexin Kong. 2023. "Two-Stage Optimization of Airport Ferry Service Delay Considering Flight Uncertainty." *European Journal of Operational Research* 307(3):1103–16. doi: 10.1016/j.ejor.2022.09.023.

He, Fang, Yafeng Yin, and Jing Zhou. 2015. "Deploying Public Charging Stations for Electric Vehicles on Urban Road Networks." *Transportation Research. Part C, Emerging Technologies* 60:227–40. doi: 10.1016/j.trc.2015.08.018.

Hiermann, Gerhard, Jakob Puchinger, Stefan Ropke, and Richard F. Hartl. 2016. "The Electric Fleet Size and Mix Vehicle Routing Problem with Time Windows and Recharging Stations." *European Journal of Operational Research* 252(3):995–1018. doi: 10.1016/j.ejor.2016.01.038.

Hof, Julian, Michael Schneider, and Dominik Goeke. 2017. "Solving the Battery Swap Station Location-Routing Problem with Capacitated Electric Vehicles Using an AVNS Algorithm for Vehicle-Routing Problems with Intermediate Stops." *Transportation Research Part B: Methodological* 97:102–12. doi: 10.1016/j.trb.2016.11.009.

Hõimoja, H., A. Rufer, G. Dziechciaruk, and A. Vezzini. 2012. "An Ultrafast EV Charging Station Demonstrator." Pp. 1390–95 in *Automation and Motion International Symposium on Power Electronics Power Electronics, Electrical Drives*.

Kang, D. , Hu, F. , & Levin, M. W. . 2022. "Impact of automated vehicles on traffic assignment, mode split, and parking behavior." *Transportation research, Part D. Transport and environment*(Mar.), 104. doi: 10.1016/j.trd.2022.103200.

Kang, D., & Levin, M. W. 2021. "Maximum-stability dispatch policy for shared autonomous vehicles." *Transportation Research Part B: Methodological*, 148, 132-151. doi: 10.1016/j.trb.2021.04.011.

Kang, D., Li, Z. and Levin, M. W. 2022. "Evasion planning for autonomous intersection control based on an optimized conflict point control formulation." *Journal of Transportation Safety & Security*, 14(12), pp. 2074–2110. doi: 10.1080/19439962.2021.1998939.

Khalkhali, Hassan, and Seyed Hossein Hosseinian. 2020. "Multi-Class EV Charging and Performance-Based Regulation Service in a Residential Smart Parking Lot." *Sustainable Energy, Grids and Networks* 22:100354-. doi: 10.1016/j.segan.2020.100354.

Lech, Norbert, and Piotr Nikończuk. 2022. "The Method of Route Optimization of Electric Vehicle." *Procedia Computer Science* 207:4454–62. doi: 10.1016/j.procs.2022.09.509.

Levin, Michael, and Di Kang. 2023. "A Multiclass Link Transmission Model for Dynamic Network Loading of Mixed Legacy and Automated Vehicle Flow." *Journal of Transportation Engineering Part A Systems 149(10).* doi:10.1061/JTEPBS.TEENG-7940.

Li, Cheng, Tao Ding, Xiyuan Liu, and Can Huang. 2018. "An Electric Vehicle Routing Optimization Model With Hybrid Plug-In and Wireless Charging Systems." *IEEE Access* 6(na):27569–78. doi: 10.1109/ACCESS.2018.2832187.

Lin, Bo, Bissan Ghaddar, and Jatin Nathwani. 2022. "Deep Reinforcement Learning for the Electric Vehicle Routing Problem With Time Windows." *IEEE Transactions on Intelligent Transportation Systems* 23(8):11528–38. doi: 10.1109/TITS.2021.3105232.

Ma, Changxi, Wei Hao, Ruichun He, Xiaoyan Jia, Fuquan Pan, Jing Fan, and Ruiqi Xiong. 2018. "Distribution Path Robust Optimization of Electric Vehicle with Multiple Distribution Centers." *PLOS ONE* 13(3):e0193789. doi: 10.1371/journal.pone.0193789.

Macrina, Giusy, Gilbert Laporte, Francesca Guerriero, and Luigi Di Puglia Pugliese. 2019. "An Energy-Efficient Green-Vehicle Routing Problem with Mixed Vehicle Fleet, Partial Battery Recharging and Time Windows." *European Journal of Operational Research* 276(3):971–82. doi: 10.1016/j.ejor.2019.01.067.

Montoya, Alejandro, Christelle Guéret, Jorge E. Mendoza, and Juan G. Villegas. 2017. "The Electric Vehicle Routing Problem with Nonlinear Charging Function." *Transportation Research Part B: Methodological* 103:87–110. doi: 10.1016/j.trb.2017.02.004.

Tadayon-Roody, Pooya, Maryam Ramezani, and Hamid Falaghi. 2021. "Multi-objective Locating of Electric Vehicle Charging Stations Considering Travel Comfort in Urban Transportation System." *IET Generation, Transmission & Distribution* 15(5):960–71. doi: 10.1049/gtd2.12072.

Tang, Qinghua, Demin Li, Yihong Zhang, and Xuemin Chen. 2023. "Dynamic Path-Planning and Charging Optimization for Autonomous Electric Vehicles in Transportation Networks." *Applied Sciences* 13(9):5476-. doi: 10.3390/app13095476.

Tu, Wei, Qingquan Li, Zhixiang Fang, Shih-lung Shaw, Baoding Zhou, and Xiaomeng Chang. 2016. “Optimizing the Locations of Electric Taxi Charging Stations: A Spatial–Temporal Demand Coverage Approach.” *Transportation Research. Part C, Emerging Technologies* 65:172–89. doi: 10.1016/j.trc.2015.10.004.

Turkensteen, Marcel. 2017. “The Accuracy of Carbon Emission and Fuel Consumption Computations in Green Vehicle Routing.” *European Journal of Operational Research* 262(2):647–59. doi: 10.1016/j.ejor.2017.04.005.

Uhrig, M., L. Weiss, M. Suriyah, and T. Leibfried. 2015. “EVS 28 KINTEX , Korea , May 3-6 , 2015 E-Mobility in Car Parks – Guidelines for Charging Infrastructure Expansion Planning and Operation Based on Stochastic Simulations.”

Van Baaren, Edzard. 2019. “The Feasibility of a Fully Electric Aircraft Towing System.”

Wang, Ying-Wei, and Chuah-Chih Lin. 2013. “Locating Multiple Types of Recharging Stations for Battery-Powered Electric Vehicle Transport.” *Transportation Research. Part E, Logistics and Transportation Review* 58:76–87. doi: 10.1016/j.tre.2013.07.003.

Wen, M., E. Linde, S. Ropke, P. Mirchandani, and A. Larsen. 2016. “An Adaptive Large Neighborhood Search Heuristic for the Electric Vehicle Scheduling Problem.” *Computers & Operations Research* 76:73–83. doi: 10.1016/j.cor.2016.06.013.

Wu, Fei, and Ramteen Sioshansi. 2017. “A Stochastic Flow-Capturing Model to Optimize the Location of Fast-Charging Stations with Uncertain Electric Vehicle Flows.” *Transportation Research. Part D, Transport and Environment* 53(C):354–76. doi: 10.1016/j.trd.2017.04.035.

Yao, Enjian, Tong Liu, Tianwei Lu, and Yang Yang. 2020. “Optimization of Electric Vehicle Scheduling with Multiple Vehicle Types in Public Transport.” *Sustainable Cities and Society* 52:101862-. doi: 10.1016/j.scs.2019.101862.

Zhang, Ying, Baltabay Aliya, Yutong Zhou, Ilsun You, Xin Zhang, Giovanni Pau, and Mario Collotta. 2018. “Shortest Feasible Paths with Partial Charging for Battery-Powered Electric Vehicles in Smart Cities.” *Pervasive and Mobile Computing* 50:82–93. doi: 10.1016/j.pmcj.2018.08.001.

Zhao, P. X., W. Q. Gao, X. Han, and W. H. Luo. 2019. “Bi-Objective Collaborative Scheduling Optimization of Airport Ferry Vehicle and Tractor.” *International Journal of Simulation Modelling* 18(2):355–65. doi: 10.2507/IJSIMM18(2)CO9.

Zhu, Shurui, Huijun Sun, and Xin Guo. 2022. “Cooperative Scheduling Optimization for Ground-Handling Vehicles by Considering Flights’ Uncertainty.” *Computers & Industrial Engineering* 169:108092. doi: 10.1016/j.cie.2022.108092.

Zweistra, Marisca, Stan Janssen, and Frank Geerts. 2020. “Large Scale Smart Charging of Electric Vehicles in Practice.” *Energies (Basel)* 13(2):298-. doi: 10.3390/en13020298.